\documentclass[10pt, reqno]{amsart} 
\usepackage{a4wide}
\usepackage{amssymb}
\usepackage{amsmath}
\usepackage{amsthm}
\usepackage{tikz}
\usepackage{amscd}
\usepackage{tikz-cd}
\usepackage[curve]{xypic} \usepackage{hyperref} \usepackage[all]{xy}
\usepackage{color} \theoremstyle{plain}
\usepackage[margin=1in]{geometry}

\newcommand{\id}{\operatorname{id}}
\newcommand{\gl}{\operatorname{GL}}
\newcommand{\sch}[1]{\operatorname{{\bf  #1}}}
\newcommand{\integers}[1]{\operatorname{\mathfrak{o}_{{\it #1}}}}
\newcommand{\ideal}[1]{\operatorname{\mathfrak{p}_{{\it #1}}}}
\newcommand{\ind}{\operatorname{ind}}

\newcommand{\res}{\operatorname{res}}
\newcommand{\ho}{\operatorname{Hom}}

\newcommand{\g}[2]{\operatorname{GL}_{#1}(#2)}

\newcommand{\tr}{\operatorname{tr}}
\newcommand{\End}{\operatorname{End}}

\newcommand{\Aut}{\operatorname{Aut}}

\newtheorem{notation}{Notation}[subsection]
\newtheorem{theorem}{Theorem}[subsection]

\newtheorem{lemma}[theorem]{Lemma}
\newtheorem{remark}[theorem]{Remark}
\newtheorem{proposition}[theorem]{Proposition}

\newtheorem{definition}[theorem]{Definition}

\newtheorem{hypothesis}{Hypothesis}[section]
 
\makeatletter
\newcommand{\proofpart}[2]{\par
  \noindent\emph{Part #1: #2}\par\nobreak
  \@afterheading
}
\makeatother
\title{On classification of typical representations for ${\rm GL}_3(F)$.}
\author{\large{ Santosh Nadimpalli }}
\date{\today}
\begin{document}
{\let\thefootnote\relax\footnote{{Keywords: Representation theory of
      $p$-adic groups. Bushnell--Kutzko theory. Typical
      representations. }}}
\begin{abstract}
  Let $F$ be any non-Archimedean local field with residue field of
  cardinality $q_F$.  In this article, we obtain a classification of
  typical representations for the Bernstein components associated to
  the inertial classes of the form
  $[\g{n}{F}\times F^\times, \sigma\otimes\chi]$ with $q_F>2$, and for
  the principal series components with $q_F>3$. With this we complete
  the classification of typical representations for $\g{3}{F}$, for
  $q_F>2$.
\end{abstract}
\maketitle
\section{Introduction}
Let $F$ be any non-Archimedean local field with residue field $k_F$ of
cardinality $q_F$.  Let $\mathcal{A}_n$ be the set of isomorphism
classes of irreducible smooth complex representations of
$\g{n}{F}$. The theory of Bernstein decomposition gives a natural
partition of the set $\mathcal{A}_n$
$$\mathcal{A}_n=\prod_{s\in \mathcal{B}_n}\mathcal{A}_n(s).$$
Here, the set $\mathcal{A}_n(s)$ is defined in terms of parabolic
induction.  The parameter $s$ is the inertial class containing the
cuspidal support of an irreducible smooth representation of $\g{n}{F}$
(see Section \ref{section_prelimi}). In the context of the local
Langlands correspondence for $\g{n}{F}$, the parameter $s$ determines
the isomorphism class of the restriction to the inertia subgroup $I_F$
of the Weil--Deligne representation associated by the local Langlands
correspondence.

The reciprocity map of the local class field theory gives an
isomorphism between the abelianization of $I_F$ and
$\mathfrak{o}_F^{\times}$, the group of units of the ring of integers
of $F$. It is natural to ask for a relation between the
representations of $I_F$, which can be extended to a Weil--Deligne
representation, and the representations of the maximal compact
subgroup $\g{n}{\integers{F}}$. One natural way would be to understand
the cuspidal support of a smooth irreducible representation from its
restriction to $\g{n}{\integers{F}}$. Indeed, in several arithmetic
applications (see \cite{Henniart-gl_2}, \cite{geo_breuil_mez}) it is
desired to construct irreducible smooth representations $\tau_s$ of
the maximal compact subgroup $\g{n}{\integers{F}}$ such that, for any
irreducible smooth representation $\pi$ of $\g{n}{F}$,
$$\ho_{ \g{n}{\integers{F}}}(\tau_s, \pi)\neq 0 \Rightarrow \pi
\in \mathcal{A}_n(s).$$
Such a representation $\tau_s$ is called a {\it typical
  representation} for $s$.

The existence of typical representations, for any $s$, follows from the
theory of types developed by Bushnell and Kutzko. For all
$s\in \mathcal{B}_n$, Bushnell and Kutzko explicitly constructed pairs
of the form $(J_s, \lambda_s)$, where $J_s$ is a compact open subgroup
of $\g{n}{F}$ and $\lambda_s$ is an irreducible smooth representation
of $J_s$ such that, for any irreducible smooth representation $\pi$ of
$\g{n}{F}$, we have
$$\ho_{ J_s}(\lambda_s, \pi)\neq 0 \Leftrightarrow \pi
\in \mathcal{A}_n(s).$$
We may assume that $J_s\subseteq \g{n}{\integers{F}}$. 
It follows from Frobenius reciprocity that any irreducible
sub representation of 
\begin{equation}\label{typical}
\ind_{J_s}^{\g{n}{\integers{F}}}\lambda_s
\end{equation}
is a typical representation for $s$. In general, the representation
\eqref{typical} is not irreducible; therefore, we cannot expect to have
a uniqueness result on typical representations for a general $s$. Now,
it is natural to ask whether there exist any other typical
representations which do not occur as subrepresentations of
\eqref{typical}. Hence, the question we are interested in is the
classification of all typical representations. In this article, we
achieve this classification for certain non-cuspidal inertial classes
of $\g{n}{F}$--including all non-cuspidal inertial classes of
$\g{3}{F}$--which give the classification of the typical
representations for all inertial classes of $\g{3}{F}$, when $q_F>2$.

Henniart (see \cite{Henniart-gl_2}) classified typical representations
for all inertial classes of $\g{2}{F}$.  Later Pa\v{s}k\={u}nas (see
\cite{Paskunas-uniqueness}) classified typical representations
occurring in the cuspidal representations of $\g{n}{F}$, for
$n\geq 3$. It turns out that there exists a unique typical
representation occurring in each cuspidal representation. Typical
representations for depth-zero inertial classes of $\g{n}{F}$ are
classified by the author in the article
\cite{level_zero_gl_n_types}. We refer to the articles of Latham
\cite{latham_sl_2}, \cite{Latham2018}, and \cite{latham_2} on typical
representations for cuspidal representations of ${\rm SL}_2(F)$ (the
tame case), cuspidal representations of ${\rm SL}_n(F)$ (the tame
case) and depth-zero cuspidal representations respectively. We also
refer to the article \cite{nevins_latham} for some results on the
typical representations for the toral cuspidal representations. The
classification of the typical representations for the non-cuspidal
inertial classes remains an open question even for ${\rm GL}_n(F)$
in the higher depth case. In this article we prove the following results.
\begin{theorem} 
  Let $n>2$ be an integer and $q_F>2$. Let $s$ be an inertial class of
  the form $[\g{n-1}{F}\times F^\times, \sigma\otimes \eta]$, where
  $\sigma$ is a cuspidal representation of $\g{n-1}{F}$ and $\eta$ is
  a character of $F^\times$.  Any typical representation $\tau_s$ for
  $s$ is isomorphic to $\ind_{J_s}^{\g{n}{\integers{F}}}\lambda_s$,
  where $(J_s, \lambda_s)$ is a Bushnell--Kutzko type for $s$.
\end{theorem}
\begin{theorem} 
  Let $n>2$ be an integer, and let $q_F>3$ if $n\neq 3$ and let
  $q_F>2$ if $n=3$. Let $s$ be an inertial class of the form $[T,
  \chi]$, where $T$ is a maximal $F$-split torus contained in
  $\g{n}{F}$ and $\chi$ is a smooth character of $T$.  Any typical
  representation $\tau_s$ for $s$ is a subrepresentation of
  $\ind_{J_s}^{\g{n}{\integers{F}}}\lambda_s$, where $(J_s,
  \lambda_s)$ is a Bushnell--Kutzko type for $s$.
\end{theorem}
Any non-cuspidal inertial class of $\g{3}{F}$ is of the above
form. Combined with the result of Pa\v{s}k\={u}nas on the unicity of
typical representation for cuspidal inertial classes, we prove the
following theorem.
\begin{theorem}
  Let $q_F>2$. Let $s$ be any  inertial class of
  $\g{3}{F}$. Any typical representation $\tau_s$ for $s$ occurs as a
  subrepresentation of $\ind_{J_s}^{\g{3}{\integers{F}}}(\lambda_s)$,
  where $(J_s, \lambda_s)$ is a Bushnell--Kutzko type for $s$.
\end{theorem}
In our analysis we will also obtain a certain multiplicity result on
the typical representations $\tau_s$.

We briefly explain the method of proof. Let $M$ be a Levi subgroup of
an $F$-parabolic subgroup $P$ of $\g{n}{F}$.  Let $\sigma$ be a
cuspidal representation of $M$. Let $\tau$ be the unique
$M\cap \g{n}{\integers{F}}$ typical representation contained in
$\sigma$. The uniqueness of $\tau$ is a result of Pa\v{s}k\={u}nas in
the article \cite{Paskunas-uniqueness}.  In order to classify typical
representations, we begin by decomposing the representation
$$\res_{\g{n}{\integers{F}}}i_{P}^{\g{n}{F}}(\sigma).$$
It follows from the results of the article
\cite{level_zero_gl_n_types} that 
$$\res_{\g{n}{\integers{F}}}i_{P}^{\g{n}{F}}(\sigma)=\ind_{\g{n}{\integers{F}}\cap
  P}^{\g{n}{\integers{F}}}\tau\oplus \Gamma,$$ where any irreducible
$\g{n}{\integers{F}}$-subrepresentation of $\Gamma$ is not a typical
representation. 

We then construct compact open subgroups of $\g{n}{\integers{F}}$,
denoted by $H_m$, for $m\geq 1$ such that
$$H_{m+1}\subset H_m,\ \text{for all} \ m\geq 1.$$ 
and 
$$\bigcap_{m\geq 1} H_m=P\cap \g{n}{\integers{F}}.$$
We will also show that $\tau$ extends as a representation of $H_m$,
for $m\geq 1$. We will then show that any $\g{n}{\integers{F}}$
subrepresentation of the representation
$\ind_{H_{m+1}}^{\g{n}{\integers{F}}}\tau/\ind_{H_m}^{\g{n}{\integers{F}}}\tau$
is not typical. The group $H_1$ will be close enough to the compact
subgroup $J_s$ in a Bushnell--Kutzko type $(J_s, \lambda_s)$. With
some more additional work, similar to the above procedure, we complete
the classification of typical representations. This requires the
analysis of the induced representation $\ind_{H_{m+1}}^{H_m}\id$. We
will also require some subtle aspects in the theory of
Bushnell--Kutzko types. In fact, the monumental theory of
Bushnell--Kutzko is the fundamental basis for this article.
\subsection{Acknowledgements.}
This article is based on chapter $4$ and $5$ of my Orsay thesis. I
would like to thank my thesis advisor Guy Henniart for suggesting this
problem and numerous discussions. I thank Corinne Blondel for pointing
out several corrections and improvements in my thesis. I want to thank
Shaun Stevens for his interest in this work. I thank the anonymous
referee for helpful suggestions.
\section{Preliminaries}\label{section_prelimi}
For any ring $R$ with unity, we denote by ${\rm Mat}_{n\times m}(R)$
the set of $n\times m$ matrices with entries in $R$. For any matrix
$X$, we denote by $X^{\text{T}}$ the transpose of $X$. The identity
matrix in ${\rm Mat}_{n\times n}(R)$ is denoted by $\id_n$ or by
$1_n$.

Let $F$ be any non-Archimedean local field with its ring of integers
$\mathfrak{o}_F$. Let $\mathfrak{p}_F$ be the maximal ideal of
$\mathfrak{o}_F$, and $\varpi_F$ be an uniformiser of $F$. Let $k_F$
be the residue field of $F$. We denote by $q_F$ the cardinality of
$k_F$. For any character $\chi$ of $F^{\times}$, we denote by
$l(\chi)$ the level of $\chi$, i.e., the least positive integer $m$
such that $1+\ideal{F}^m$ is contained in the kernel of $\chi$. Note
that the level $l(\chi)$ of an unramified character $\chi$ is still
$1$. We denote by $\nu_F:F^\times \rightarrow\mathbb{Z}$ the
normalised valuation of $F$.

All representations in this article are defined over  complex vector
spaces. 

Let $G$ be the $F$-rational points of a connected reductive algebraic
group $\sch{G}$ defined over $F$. Let $\mathcal{R}(G)$ be the category
of smooth representations of $G$. For any closed subgroup $H$ of $G$,
we denote by $\ind_H^G$ the compact induction functor from
$\mathcal{R}(H)$ to $\mathcal{R}(G)$. Let $P$ be the group of
$F$-rational points of any $F$-parabolic subgroup of $\sch{G}$. Let
$M$ be a Levi subgroup of $P$. We denote by $i_{P}^G$ the normalised
parabolic induction functor from $\mathcal{R}(M)$ to $\mathcal{R}(G)$.

Let $H_1$ and $H_2$ be two groups and $\tau_1$ and
$\tau_2$ be any representations of $H_1$ and $H_2$ respectively. We
denote by $\tau_1\boxtimes\tau_2$ the tensor product representation of
$H_1\times H_2$. If $H_1=H_2$, then the representation
$\tau_1\otimes\tau_2$ is the tensor product representation of $H_1$.

For any positive integer $n$, the group ${\rm GL}_n(F)$ is denoted by
$G_n$ and the group ${\rm GL}_n(\mathfrak{o}_F)$ is denoted by
$K_n$. The principal congruence subgroup of $K_n$ of level $m$ is
denoted by $K_n(m)$, for $m\geq 1$. Let $I$ be a sequence of positive
integers $(n_1, n_2, \dots, n_r)$ such that $n_1+n_2+\dots+n_r=n$. For
any ring $R$ with unity, we denote by $P_I(R)$ (resp. $\bar{P}_I(R)$)
the group of invertible block upper (resp. lower) triangular matrices
of the type $I$. Let $U_I(R)$ (resp. $\bar{U}_I(R)$) be the group of
block upper (lower) unipotent matrices of the type $I$. Let $M_I(R)$
be the block diagonal matrices of the type $I$. If $I=(1,1,\dots,1)$,
the groups $P_I(R)$, $M_I(R)$ and $U_I(R)$ and $\bar{U}_{I}(R)$ are
denoted by $B_n(R)$, $T_n(R)$ and $U_n(R)$ and $\bar{U}_n(R)$
respectively.  When $R=F$, we drop the symbol $R$, i.e., $P_I(R)$ will
be denoted by $P_I$ etc. We have $P_I=M_IU_I$ and
$\bar{P}_I=M_I\bar{U}_I$.

As an example, when $I=(n-1, 1)$, the group $P_{I}$ in the block form
is given by:
$$\begin{pmatrix}G_{n-1}&{\rm Mat}_{n-1\times 1}(F)\\
  0&F^\times\end{pmatrix}.$$ In the
block form the groups $M_I$ and $U_I$ are given by
$$M_I=\begin{pmatrix}G_{n-1}&0\\
  0&F^\times\end{pmatrix}\ \text{and}\  U_I=
\begin{pmatrix}1_{n-1}&{\rm Mat}_{n-1\times 1}(F)\\
    0&1\end{pmatrix}.$$
  We identify the
group $M_{(n-1, 1)}$ with the group $G_{n-1}\times G_1$. Any
irreducible smooth representation of $M_{(n-1,1)}$ is identified with
$\sigma\boxtimes\chi$, where $\sigma$ is an irreducible smooth
representation of $G_{n-1}$ and $\chi$ is a character of $G_1$. 

We briefly recall the theory of Bernstein decomposition. Let $B(G)$ be
the set of pairs $(M, \sigma)$, where $M$ is a Levi subgroup of an
$F$-parabolic subgroup $P$ of $G$, and $\sigma$ is an irreducible
cuspidal representation of $M$. The pairs $(M_1, \sigma_1)$ and
$(M_2, \sigma_2)$ in $B(G)$ are said to be inertially equivalent 
if and only if there exist an element $g\in G$ and an unramified
character $\chi$ of $M_2$ such that 
$$M_1=gM_2g^{-1}\ \text{and}\ \sigma_1^g\simeq \sigma_2\otimes\chi.$$
We denote by $\mathcal{B}_G$ the set of equivalence classes, called 
{\it inertial classes}. 

Any irreducible smooth representation $\pi$ of $G$ occurs as a
sub-representation of a parabolic induction $i_{P}^{G}(\sigma)$, where
$\sigma$ is an irreducible cuspidal representation of a Levi subgroup
$M$ of $P$. The pair $(M, \sigma)$ is well determined up to
$G$-conjugation. We call the class $s=[M, \sigma]$ the {\it inertial
  support} of $\pi$. We denote by $\mathcal{I}(\pi)$ the inertial
support of $\pi$. For any inertial class $s=[M, \sigma]$, we denote by
$\mathcal{R}_s(G)$ the full sub-category of $\mathcal{R}(G)$
consisting of smooth representations all of whose irreducible
sub-quotients have inertial support $s$.  It is shown by Bernstein in
\cite{le_centre_bernstein} that the category $\mathcal{R}(G)$
decomposes as a direct product of $\mathcal{R}_s(G)$. The category
$\mathcal{R}_s(G)$ is called a {\it Bernstein component} associated to
the inertial class $s$.  In particular, every smooth representation can
be written as a direct sum of objects in the categories
$\mathcal{R}_s(G)$. We denote by $\mathcal{A}_n(s)$ the set of
isomorphism classes of irreducible representations in
$\mathcal{R}_s(G_n)$.
\begin{definition}\label{definition_main}
  Let $s$ be an inertial class for $G_n$. An irreducible smooth
  representation $\tau$ of $K_n$ is called a typical representation
  for $s$, if for any irreducible smooth representation $\pi$ of $G$,
  we have
  $$\ho_{K_n}(\tau, \pi)\neq 0 \implies
  \mathcal{I}(\pi)=s.$$ A non typical representation is called an
  atypical representation.
\end{definition}
For any inertial class $s$ of $G_n$, the existence of a typical
representation can be deduced from the theory of types developed by
Bushnell and Kutzko in the book \cite{Orrangebook} and the article
\cite{Bushnell-kutzko-Semisimpletypes}. Bushnell and Kutzko
constructed explicit pairs (called types) $(J_s, \lambda_s)$, where
$J_s$ is a compact open subgroup of $\g{n}{F}$, and $\lambda_s$ is an
irreducible smooth representation of $J_s$. The pair $(J_s,\lambda_s)$
satisfies the condition that, for any irreducible smooth
representation $\pi$ of $G$, we have
$$\ho_{J_s}(\pi, \lambda_s)\neq 0\ \Leftrightarrow\ \mathcal{I}(\pi)=s.$$ 
The group $J_s$ can be arranged to be a
subgroup of $\g{n}{\integers{F}}$ by conjugating with an element of
$\g{n}{F}$. Hence we assume that $J_s\subseteq
\g{n}{\integers{F}}$.
It follows from Frobenius reciprocity that any irreducible
sub-representation of
\begin{equation}\label{bushnell_max_induction}
\ind_{J_s}^{\g{n}{\integers{F}}}(\lambda_s)
\end{equation}
is a typical representation. The irreducible sub representations of
(\ref{bushnell_max_induction}) are classified by Schneider and Zink in
\cite[Section 6, $T_{K,\lambda}$ functor]{SchneiderKtypes}.

For $s=[G_n, \sigma]$, Pa\v{s}k\={u}nas in the article
\cite[Theorem 8.1]{Paskunas-uniqueness} showed that up to isomorphism
there exists a unique typical representation for $s$. More precisely,
\begin{theorem}[Pa\v{s}k\={u}nas]\label{paskunas_prelim_rev}
  Let $n\geq 1$ be an integer and $\sigma$ be an irreducible cuspidal
  representation of $G_n$. Let $(J_s, \lambda_s)$ be a Bushnell-Kutzko
  type for the inertial class $s=[G_n, \sigma]$ with
  $J_s\subseteq K_n$. The representation
$$\ind_{J_s}^{K_n}(\lambda_s)$$
is the unique typical representation for the inertial class $[G_n,
\sigma]$. The representation $\ind_{J_s}^{K_n}(\lambda_s)$ occurs with
a multiplicity one in $\sigma$.
\end{theorem}
In this article, we classify typical representations for ${\rm
  GL}_3(F)$ in terms of Bushnell--Kutzko types. We first obtain a
classification of typical representations for the inertial classes
$[M_{(n-1, 1)}, \sigma\boxtimes\eta]$ and $[T_n, \chi]$, where $\eta$
and $\chi$ are characters of $F^\times$ and $T_n$ respectively. We
will use some basic results from the article
\cite{level_zero_gl_n_types} and we recall some of these results.
\begin{lemma}\label{lemma_prelim_twisting}
  Let $\chi$ be a character of $G_n$ and let $\tau$ be a typical
  representation for an inertial class $s=[M, \sigma]$ of $G_n$.  The
  representation $\tau\otimes\chi$ is a typical representation for the
  inertial class $[M, \sigma\otimes\chi]$.
\end{lemma}
\begin{proof}
  We refer to \cite[Lemma 2.7]{level_zero_gl_n_types} for a proof.
\end{proof}
Let $P$ be any parabolic subgroup of $G_n$ with a Levi subgroup $M$
and $U$ be the unipotent radical of $P$. Let $\bar{U}$ be the
unipotent radical of the opposite parabolic subgroup of $P$ with
respect to $M$. Let $J_1$ and $J_2$ be two compact open subgroups of
$K_n$ such that $J_1$ contains $J_2$.  Suppose $J_1$ and $J_2$ both
satisfy the Iwahori decomposition with respect to $P$ and $M$. With
$J_1\cap U=J_2\cap U$ and $J_1\cap \bar{U}=J_2\cap \bar{U}$. Let
$\lambda$ be an irreducible smooth representation of $J_2$ which
admits an Iwahori decomposition i.e. $J_2\cap U$ and $J_2\cap \bar{U}$
are contained in the kernel of $\lambda$.
\begin{lemma}\label{prelim_transitivity_of_covers}
  The representation $\ind_{J_2}^{J_1}(\lambda)$ is the extension of
  the representation $\ind_{J_2\cap M}^{J_1\cap M}(\lambda)$ such that
  $J_1\cap U$ and $J_1\cap \bar{U}$ are contained in the kernel of the
  extension.
\end{lemma}
\begin{proof}
  The lemma is well known and frequently used when dealing with the
  formalism of $G$-covers. We refer to \cite[Lemma
  2.6]{level_zero_gl_n_types} for a proof. 
\end{proof}
Let $t_i=[M_i, \Theta_i]$ be an inertial class of $G_{n_i}$,
for $1\leq i\leq r$.  Let $\sigma_i$ be a smooth representation from
$\mathcal{R}_{t_i}({G_{n_i}})$. We suppose
$$\res_{K_{n_i}}\sigma_i=\tau_i^0\oplus\tau_i^1,$$ for
$1\leq i\leq r$, such that irreducible $K_{n_i}$-subrepresentations of
$\tau_i^1$ are atypical.  We denote by $t$ the inertial class
$$[M_1\times M_2 \times \dots \times M_r,
\Theta_1\boxtimes\Theta_2\boxtimes\dots\boxtimes\Theta_r]$$ of
$G_n$. 
The inertial class $t$ is independent of the choice of representatives 
$(M_i, \Theta_i)$. Let $\tau_I^0=\boxtimes_{i=1}^{r}\tau_i^0$ and 
$\sigma_I=\boxtimes_{i=1}^{r}(\sigma_i)$.
\begin{lemma}\label{lemma_prelim_prelim_elim}
  The representation $\ind_{P_I\cap K_n}^{K_n}(\tau_I^0)$ admits a
  complement in $\res_{K_n}i_{P_I}^{G_n}(\sigma_I)$ with all its
  irreducible sub-representations atypical.
\end{lemma}
\begin{proof}
  We refer to \cite[Proposition 2.3]{level_zero_gl_n_types} for a
  proof.
\end{proof}
In particular, if $\Theta_i=\sigma_i$ is a cuspidal representation,
then from Theorem \ref{paskunas_prelim_rev} we have
$\res_{K_{n_i}}\sigma_i=\tau_i^0\oplus \tau_i^1$, where $\tau^0_i$ is
the unique typical representation for the inertial class $[G_{n_i},
\sigma_i]$. Hence, any typical representation for the inertial class
$t$ occurs as a sub-representation of $\ind_{P_I\cap
  K_n}^{K_n}\sigma_I$.
\begin{lemma}\label{standard_form_rev}
  Let $s=[M, \sigma]$ be any inertial class of $G_n$. Then there
  exists a partition $I$ of $n$ and a cuspidal representation
  $\sigma_I$ of $M_I$ such that $s=[M_I, \sigma_I]$.
\end{lemma}
\begin{proof}
  We refer to \cite[Section 2.2, page no. 5]{level_zero_gl_n_types}
  for a proof.
\end{proof}
The following result is useful in understanding some stabilisers in
the later part of this article.  The space ${\rm Mat}_{n\times
  m}(k_F)$ is equipped with an action of
$M_{(m,n)}(k_F)=\g{m}{k_F}\times \g{n}{k_F}$ given by $(g_1,
g_2)U=g_2Ug_1^{-1}$, for $U\in{\rm Mat}_{n\times m}(k_F)$. We also
have a $M_{(m,n)}(k_F)$ action on the set of matrices ${\rm
  Mat}_{m\times n}(k_F)$ by setting $(g_1,g_2)V=g_1Vg_2^{-1}$, for
$V\in{\rm Mat}_{n\times m}(k_F)$. Let $\psi$ be a non-trivial
character of the additive group $k_F$. We define a pairing $B$ between
${\rm Mat}_{m\times n}(k_F)$ and ${\rm Mat}_{n\times m}(k_F)$ by
defining $B(V,U)=\psi\circ\tr(VU)$. Let $T$ be the map from ${\rm
  Mat}_{m\times n}(k_F)$ and ${\rm Mat}_{n\times m}(k_F)^{\wedge}$
defined by
$$T(V)(U)=B(V,U).$$
\begin{lemma}\label{prelim_duality_isom}
  The map $T$ is an $M_{(m,n)}(k_F)$-equivariant isomorphism.
\end{lemma}
Let $s$ be any depth-zero inertial class $[M_I, \sigma_I]$ of
$G_n$. The group $K_n\cap M_I$ acts on the space 
$$\sigma_I^{K_n(1)\cap M_I},$$
and we denote this representation of $K_n\cap M_I$ by $\tau_I$. The
pair $(K_n\cap M_I, \tau_I)$ is a Bushnell--Kutzko type for the
inertial class $[M_I, \sigma_I]$ of $M_I$. Let $P_I(1)$ be the group
$K_n(1)(P_I\cap K_n)$ and observe that $P_I(1)\cap M_I$ is equal to
$K_n\cap M_I$. The representation $\tau_I$ extends as a representation
of $P_I(1)$ such that $P_I(1)\cap U_I$ and $P_I(1)\cap \bar{U}_I$ are
contained in the kernel of this extension.
\begin{theorem}\label{thr_depth_zero}
  Let $s=[M_I, \sigma_I]$ be any depth-zero inertial class of
  $G_n$. Any typical representation $\tau_s$ for $s$ occurs as a
  subrepresentation of $\ind_{P_I(1)}^{K_n}\tau_I$. Moreover, we have
  $$\dim_{\mathbb{C}}\ho_{K_n}(\tau_s, i_{P_I}^{G_n}\sigma_I)=
  \dim_{\mathbb{C}}\ho_{K_n}(\tau_s, \ind_{P_I(1)}^{K_n}\tau_I).$$
\end{theorem}
\begin{proof}
  We refer to \cite[Theorem 3.2]{level_zero_gl_n_types} for the proof.
  \end{proof}
\section{The inertial class with 
Levi subgroup of type \texorpdfstring{$(n-1,1)$}{}}
Let $n>1$ be any positive integer. In this section we assume that
  $I=(n-1, 1)$. Let $V$ and $V_1$ be two $F$-vector spaces of
dimensions $n-1$ and $1$ respectively. Let $P$ be the parabolic
subgroup of $GL(V\oplus V_1)$ fixing the flag
$V\subseteq V\oplus V_1$. We denote by $M$ its Levi subgroup fixing
the decomposition $V\oplus V_1$.  Hence, we have
$M=GL(V)\times GL(V_1)$. In this section, we are interested in the
classification of typical representations for inertial classes
$[M, \sigma\boxtimes\chi]$, where $\sigma$ is a cuspidal
representation of $GL(V)$, and $\chi$ is a character of $GL(V_1)$. We
will use the language of the book \cite{Orrangebook} freely in this
section. Let $(J(\mathfrak{A}, \beta), \lambda)$ be a maximal simple
(Bushnell--Kutzko) type contained in the representation $\sigma$. We
recall certain important features of this type for our purpose.
\subsection{Bushnell-Kutzko semi-simple type}\label{n_1-type}
We denote by $A$ the algebra $\End_{F}(V)$. Let
$[\mathfrak{A},l,0,\beta]$ be a simple stratum in $A$ defining the
maximal simple type $(J(\mathfrak{A}, \beta), \lambda)$. We denote by
$B$ the commutant of $E=F[\beta]$ in $A$. Let
$\mathfrak{B}=\mathfrak{A}\cap B$. Let $\mathfrak{P}$ and
$\mathfrak{D}$ be the radicals of $\mathfrak{A}$ and $\mathfrak{B}$
respectively. Given any hereditary order $\mathfrak{A}$, we define the
filtration $U^i(\mathfrak{A})$ by setting
$$U^i(\mathfrak{A})=\id+\mathfrak{P}^i,$$
for all $i\geq 1$, and $U^0(\mathfrak{A})$ is the set of units of
$\mathfrak{A}$. The type $(J(\mathfrak{A},\beta), \lambda)$ is called
maximal if $\mathfrak{B}$ is a maximal hereditary order in $B$.

The group $J(\mathfrak{A}, \beta)$ contains $U^0(\mathfrak{B})$. There
is a normal subgroup $J^1(\mathfrak{A}, \beta)$ of
$J(\mathfrak{A}, \beta)$ such that
$J^1(\mathfrak{A}, \beta)\cap U^0(\mathfrak{B})=U^1(\mathfrak{B)}$ and
$$\dfrac{U^0(\mathfrak{B})}{U^1(\mathfrak{B})}\simeq
\dfrac{J(\mathfrak{A}, \beta)}{J^1(\mathfrak{A}, \beta)}.$$ The
group $U^0(\mathfrak{B})/U^1(\mathfrak{B})$ is a general linear group
of a vector space over a finite field. The representation $\lambda$ is
an irreducible representation which is given by a tensor product
representation $\kappa\otimes\rho$, where $\kappa$ is a representation
of $J(\mathfrak{A}, \beta)$, called a $\beta$-extension (see
\cite[Chapter 5, Definition 5.2.1]{Orrangebook}), and $\rho$ is a
cuspidal representation of $U^0(\mathfrak{B})/U^1(\mathfrak{B})$
(considered as a representation of $J(\mathfrak{A}, \beta)$ through
its quotient $J(\mathfrak{A}, \beta)/J^1(\mathfrak{A}, \beta)$). We
refer to \cite[Chapter 5]{Orrangebook} for complete details of these
constructions.  For the precise definition and description see
\cite[chapter 5, Definition 5.5.10]{Orrangebook}. For simplicity, we
will denote by $J^0$ and $J^1$ the groups $J(\mathfrak{A}, \beta)$ and
$J^1(\mathfrak{A}, \beta)$ respectively.

We fix the following notations. Let $e$ and $f$ be the ramification
index and inertial index of $E$ respectively. We fix an
$\mathfrak{o}_E$-lattice chain $\mathcal{L}$ defining the hereditary
$\mathfrak{o}_E$-order $\mathfrak{B}$. Let $\mathfrak{A}$ be the
hereditary $\mathfrak{o}_F$-order defined by the lattice chain
$\mathcal{L}$, considering $\mathcal{L}$ as an
$\mathfrak{o}_F$-lattice chain. We fix a $\mathfrak{o}_E$-basis
$(w_1,w_2,\dots w_{(n-1)/ef})$ for the lattice chain $\mathcal{L}$
(see \cite[Chapter 1, 1.1.7]{Orrangebook}) and then a
$\mathfrak{o}_F$-basis for $\mathfrak{o}_Ew_i$ for $1\leq i\leq ef$;
hence, we obtain an $F$-basis $(v_1,v_2,\dots, v_{n-1}, v_{n})$ for
the vector space $V\oplus V_1$, where $v_n\in V_1$.  In this basis, we
write all our endomorphisms as matrices of ${\rm Mat}_{n\times n}(F)$.
With this basis we have $J^0\subseteq K_{n-1}$.

We are interested in the classification of typical representations for
the inertial class $[M_I, \sigma\boxtimes\chi]$. By twisting with a
character if necessary we may (and do) assume that $\chi=\id$ (see
Lemma \ref{lemma_prelim_twisting}). Let $\tau$ be any typical
representation for the above inertial class. The representation
$\ind_{K_n}^{G_n}\tau$ is a finitely generated representation of $G_n$
and hence admits an irreducible quotient $\pi$. Using the definition
of a typical representation we see that $\pi$ occurs as a sub quotient
of $i_{P_I}^{G_n}(\sigma\chi_1\boxtimes\chi_2)$, where $\chi_1$ and
$\chi_2$ are unramified characters of $G_{n-1}$ and $G_1$
respectively.

Hence, in order to classify typical representations for the inertial class
$s=[M_I, \sigma\boxtimes\id]$, it is enough to examine which
$K_n$-irreducible sub representations of
$$\res_{K_n}i_{P_I}^{G_n}(\sigma\boxtimes\id)$$
are typical for the inertial class $s$.  Let $\tau$
be the unique typical representation contained in the representation
$\sigma$. The
representation
$$\ind_{P_I\cap K_n}^{K_n}(\tau\boxtimes\id)$$
has a complement in 
$$\res_{K_n}i_{P_I}^{G_n}(\sigma\boxtimes\chi)$$
whose irreducible sub representations are atypical (see Lemma
\ref{lemma_prelim_prelim_elim}).

Now we have to look for typical representations occurring in the
representation
$$\ind_{P_I\cap K_n}^{K_n}(\tau\boxtimes\id).$$
For this purpose, we will define compact subgroups $H_m\subseteq K_n$,
for $m\geq N_0$ (for some positive integer $N_0$), such that
\begin{hypothesis}\label{auxillary_group_cond}\leavevmode
\begin{enumerate}
\item $H_{m+1}\subseteq H_m$, for $m\geq N_0$ and
  $\bigcap_{m\geq N_0}H_m=P_I(\integers{F})$,
\item the group $H_m$ has the Iwahori decomposition with respect to
  $P_I$ and its Levi subgroup $M_I$,
\item  the representation $\tau\boxtimes\id$ admits an
  extension to $H_{N_0}$ such that $H_{N_0}\cap \bar{U}_{I}$ and
  $H_{N_0}\cap U_I$ are contained in the kernel of this extension. 
\end{enumerate}
\end{hypothesis}
For any such sequence $\{H_m, m\geq N_0\}$ as above we have:
$$\ind_{P_I\cap K_n}^{K_n}(\tau\boxtimes\id)
\simeq \bigcup_{m\geq N_0}\ind_{H_m}^{K_n}(\tau\boxtimes\id).$$ Before
we start this construction it is instructive to first examine the
Bushnell-Kutzko semi-simple type for the inertial class
$[M_I, \sigma\boxtimes\id]$.

Let us recall some standard material required from
\cite{Bushnell-kutzko-Semisimpletypes}. First, let us begin with
lattice sequences.  An $\mathfrak{o}_F$-lattice sequence in a
$F$-vector space, say $V$, is a function $\Lambda$ from $\mathbb{Z}$
to the set of $\mathfrak{o}_F$-lattices in $V$ with the following
conditions on $\Lambda$:
$$\Lambda(n+1) \subseteq \Lambda(n),\ \text{for all}\ n\in \mathbb{Z};$$
and there exists an $e(\Lambda)\in \mathbb{Z}$ such
that
$$\Lambda(n+e(\Lambda))=\ideal{F}\Lambda(n),\ \text{for all}\ n\in
\mathbb{Z}.$$ An $\mathfrak{o}_F$-lattice chain is an
$\mathfrak{o}_F$-lattice sequence with the strict inclusion between
$\Lambda(n+1)$ and $\Lambda(n)$, for all $n\in \mathbb{Z}$.  One
extends the function $\Lambda$ to the set of real numbers by setting
$$\Lambda(r)=\Lambda(-[-r]),$$
for all $r\in \mathbb{R}$. Here, $[x]$ is the greatest integer less than
or equal to $x$. 

Given two $\mathfrak{o}_F$-lattice sequences $\Lambda_1$ and
$\Lambda_2$ in the vector spaces $V_1$ and $V_2$ over $F$, Bushnell
and Kutzko defined the notion of direct sum of $\Lambda_1$ and
$\Lambda_2$. Let $e=\text{lcm}(e(\Lambda_1), e(\Lambda_2))$. The
direct sum of $\Lambda_1$ and $\Lambda_2$, denoted by $\Lambda$, is an
$\mathfrak{o}_F$-lattice sequence in the vector space $V_1\oplus V_2$
given by
$$\Lambda(er)=\Lambda_1(e_1r)\oplus \Lambda(e_2r),$$
for any $r\in \mathbb{R}$.  Given an $\mathfrak{o}_F$-lattice sequence
$\Lambda$ in a vector space $V$ one can define a filtration
$\{a_{r}(\Lambda)\ | \ r\in \mathbb{R}\}$ on the algebra $\End_F(V)$
given by the equation
$$a_{r}(\Lambda)=\{x\in \End_F(V) |\
x\Lambda(i)\subseteq \Lambda(i+r)\ \forall \ i \in \mathbb{Z}\}.$$ Let
$u_0(\Lambda)$ be the group of units in the order $a_{0}(\Lambda)$
and, for $r>0$ and $r\in \mathbb{Z}$, we set $u_r(\Lambda)$ to be
$1+a_{r}(\Lambda)$.

Let $(J_s, \lambda_s)$ be a Bushnell-Kutzko type for the inertial class
$$s=[M_I, \sigma\boxtimes\id].$$
The group $J_s$ satisfies the Iwahori decomposition with respect to the
parabolic subgroup $P_I$ and the Levi subgroup $M_I$. Let
$[\mathfrak{A}, l,0, \beta]$ be a simple stratum defining the
type $(J^0, \lambda)$ for the inertial class $[G_{n-1}, \sigma]$. The
order $\mathfrak{A}$ is defined by a lattice  chain $\Lambda_1$ with
values in sub-lattices of $\mathfrak{o}_F^n$. We denote by
$\Lambda_2$ the lattice chain defined by
$\Lambda_2(i)=\mathfrak{p}_F^i$, for all $i\in \mathbb{Z}$. We
have
\begin{enumerate}
\item $J_s\cap U_I=u_0(\Lambda_1\oplus\Lambda_2)\cap U_I$, 
\item $J_s\cap M_I=J^0\times \mathfrak{o}_F^{\times}$,
\item $J_s\cap \bar{U}_I=u_{l+1}(\Lambda_1\oplus\Lambda_2)\cap\bar{U}_I$,  
\item the restriction of $\lambda_s$ to $J_s\cap M_I$ is isomorphic to
  $\lambda\boxtimes \id$, and the groups $J_s\cap \bar{U}_I$ and
  $J_s\cap {U}_I$ are contained in the kernel of $\lambda_s$.
\end{enumerate}
We refer to \cite [Section 8, paragraph
8.3.1]{Bushnell-kutzko-Semisimpletypes} for the construction of the
above Bushnell--Kutzko type.

Now we make an explicit calculation of the groups
$u_{l+1}(\Lambda_1\oplus\Lambda_2)\cap \bar{U}_I$ and
$u_{0}(\Lambda_1\oplus\Lambda_2)\cap {U}_I$.  Note that the
period of the direct sum $\Lambda_1\oplus \Lambda_2$ is the least
common multiple of the period of the two lattice sequences
$\Lambda_1$ and $\Lambda_2$. Hence, the period of the
lattice sequence $\Lambda$ is $e$, where $e$ is the
period of the lattice chain $\Lambda_1$. Let $t$ be an integer such
that $0\leq t\leq e-1$.  Let $L_0$ be the free $\mathfrak{o}_F$ module
$\mathfrak{o}_F^{(n-1)/e}$. The lattice  $\Lambda_1(t)$ is given by
:
$$\Lambda_1(t)=(L_0\oplus L_0\oplus\dots\oplus L_0)
\oplus (\varpi_FL_0\oplus \varpi_FL_0\oplus\dots\oplus\varpi_F L_0),$$
where the $L_0$ is repeated $e-t$ times, and $\varpi_FL_0$ is repeated
$t$ times. Hence,
the lattice chain $\Lambda$ is given by
$$\Lambda(0)=\Lambda_1(0)\oplus \Lambda_2(0)=
 (L_0\oplus L_0\oplus \dots\oplus L_0)\oplus \mathfrak{o}_F$$
and 
$$\Lambda(t)=\Lambda_1(t)\oplus\Lambda_2(t/e)= 
(L_0\oplus L_0\oplus\dots\oplus L_0)\oplus (\varpi_FL_0\oplus
\varpi_FL_0\oplus\dots\oplus\varpi_F L_0)\oplus \ideal{F},$$ for
$0\leq t\leq e-1$.

We note that $u_0(\Lambda)\cap U_I=U_I(\integers{F})$.  We denote
by $\bar{\mathfrak{n}}_I$ the lower nilpotent matrices of the type
$(n-1, 1)$, i.e. the Lie algebra of $\bar{U}_I$.
We then have:
$$u_{l+1}(\Lambda)\cap \bar{U}_I=
\id+(a_{l+1}(\Lambda)\cap \bar{\mathfrak{n}}_I).$$
Let $l+1=el'+r$, where $0\leq r<e$.  Since $\Lambda$ is a lattice
chain of period $e$, we deduce that
$$a_{l+1}(\Lambda)\cap \bar{\mathfrak{n}}_I=
\varpi^{l'}_F(a_{r}(\Lambda)\cap \bar{\mathfrak{n}}_I).$$ Finally, it
remains to calculate the group $a_{r}(\Lambda)\cap
\bar{\mathfrak{n}}_I$.  We note that $a_{r}(\Lambda)\cap
\bar{\mathfrak{n}}_I$ is the following set
$$\{x\in {\rm Mat}_{n\times n}(F)\cap \bar{\mathfrak{n}}_I \ |\  
x\Lambda(i)\subseteq \Lambda(i+r)\ \forall \ i\in \mathbb{Z}\}.$$ For
$r\geq 1$, the $n^{\text{th}}$ row (in block form) of an element in
$a_r(\Lambda)\cap \bar{\mathfrak{n}}_I$ is of the form $A=[M_1,
M_2,\dots, M_{e}, 0]$, where $M_i$ is a matrix of type $1\times
(n-1)/e$, for $1\leq i\leq e$ and:
 \begin{enumerate}
 \item   $M_i\in \varpi_F^2{\rm Mat}_{1\times
     (n-1)/e}(\mathfrak{o}_F)$,
 for $i\leq r-1$,  
\item $M_i\in \varpi_F{\rm Mat}_{1\times (n-1)/e}(\mathfrak{o}_F)$,
  for $i>r-1$.
\end{enumerate}
If $r=0$ and $e>1$, then we know that
$M_{i}\in \varpi_F{\rm Mat}_{1\times (n-1)/e}(\mathfrak{o}_F)$, for
$1\leq i\leq e-1$, and
$M_e\in {\rm Mat}_{1\times (n-1)/e}(\integers{F})$. If $r=0$ and
$e=1$, then we have $A\in {\rm Mat}_{1\times n}(\integers{F})$.
This description is enough for the present purposes.
\subsection{Some auxiliary groups.}
Let $m$ be a positive integer and $P_I(m)$ be the inverse image of the
group $P_I(\integers{F}/\ideal{F}^m)$ under the mod-$\ideal{F}^m$
reduction map  
$$K_n\rightarrow {\rm GL}_n(\mathfrak{o}_F/\mathfrak{p}_F^m).$$
There exists a positive integer $N_1$ such that the principal
congruence subgroup of level $N_1$ is contained in the kernel of the
representation $\tau$. The representation $\tau\boxtimes\id$ of
$M_I(\integers{F}/\ideal{F}^{N_1})$ now extends to a representation of
$P_I(N_1)$ by inflation. We note that
$P_I(N_1)\cap \bar{U}_I$ and $P_I(N_1)\cap U_I$ are both contained in the
kernel of this extension.  Now the sequence of groups $H_m=P_I(m)$ and
the representation $\tau\boxtimes\id$, for $m\geq N_1$ satisfy the
conditions in Hypothesis \ref{auxillary_group_cond}. Hence, we get that
$$\ind_{P_I\cap K_n}^{K_n}(\tau\boxtimes\id)
\simeq \bigcup_{m\geq N_1}\ind_{P_I(m)}^{K_n}(\tau\boxtimes\id).$$
We conclude that typical representations occur as sub representations of 
$$\ind_{P_I(m)}^{K_n}(\tau\boxtimes\id),$$
for some positive integer $m\geq N_1$.  

For making Mackey decompositions easier and other reasons, it is
convenient to work with a smaller subgroup $P_I^0(m)$ of $P_I(m)$. We
begin by rewriting the representation
$$\ind_{P_I(m)}^{K_n}(\tau\boxtimes\id).$$
We also require to make $N_1$ explicit.  We recall that $K_{n-1}(m)$
is the principal congruence subgroup of level $m$ of $G_{n-1}$. The
group $J^0$ contains the group $U^{[l/2]+1}(\mathfrak{A})$ and the
representation $\lambda$ restricted to the group
$U^{[l/2]+1}(\mathfrak{A})$ is a direct sum of copies of the same character
$\psi_{\beta}$. The character $\psi_\beta$ is trivial on the group
$U^{l+1}(\mathfrak{A})$. We also recall the notation that $l+1=el'+r$,
where $0\leq r\leq e-1$. We note that
$U^{l+1}(\mathfrak{A})=\id_{n-1}+\varpi_F^{l'}\mathfrak{P}_{\mathfrak{A}}^{r}$. If
$r=0$, then $K_{n-1}(1)\subseteq \mathfrak{P}_{\mathfrak{A}}^r$. If
$r>1$, then from the formulas \cite [2.5.2]{Orrangebook} we get that
$K_{n-1}(2)\subseteq \mathfrak{P}_{\mathfrak{A}}^{r}$, for $0\leq
r<e$.  This shows that the representation $\lambda$ is trivial on
$K_{n-1}(N_s)$, where $N_s$ is given by:
\begin{notation}\label{rev_n_1_begin}
  From now we fix $N_s=[(l+1)/e]+1$ if $r=0$ and $e>1$. If $r=0$ and
  $e=1$, then $N_s=l+1$. Finally,  $N_s=[(l+1)/e]+2$ if $r\geq1$. 
\end{notation} 
Let $\pi$ be the projection map
$$P_I(\integers{F})\rightarrow M_I(\integers{F}).$$
For $m\geq N_s$, we denote by $P_I^0(m)$ the group
$K_{n}(m)\pi^{-1}(J^0\times\mathfrak{o}_F^{\times})$.  Since
$K_{n}(m)\cap P_I\subseteq \pi^{-1}(J^0\times\mathfrak{o}_F^{\times})$,
the group $P_I^0(m)$ satisfies the Iwahori decomposition with respect to
the subgroup $P_I$ and the Levi subgroup $M_I$. In particular, we have
$$P_I^0(m)=(P_I^0(m)\cap U_I)(P_I^0(m)\cap M_I)(P_I^0(m)\cap \bar{U}_I).$$
Here, $P_I^0(m)\cap U_I$ is equal to $U_I(\integers{F})$,
$P_I^0(m)\cap M_I$ is equal to
$J^0\times\mathfrak{o}_F^{\times}$, and $(P_I^0(m)\cap \bar{U}_I)$ is
equal to
$K_n(m)\cap \bar{U}_I$.

We observe that $\lambda\boxtimes\id$ extends as a representation of
$P^0(m)$, for all $m\geq N_s$; the groups $P^0(m)\cap U_I$ and
$P^0(m)\cap \bar{U}_I$ are contained in the kernel of this extension.
Now the representation $\tau\boxtimes\id$ of
$K_{n-1}\times\mathfrak{o}_F^{\times}$ is isomorphic to
$$\{\ind_{J^0}^{K_{n-1}}(\lambda)\}\boxtimes\id.$$ 
Hence, we get that
$$\ind_{P_I(m)}^{K_n}
(\tau\boxtimes\id)\simeq \ind_{P_I^0(m)}^{K_n}(\lambda\boxtimes\id),$$
for all $m\geq N_s$ (we apply Lemma
\ref{prelim_transitivity_of_covers} to the groups $J_1=P_I(m)$ and
$J_2=P_I^0(m)$ and $\lambda=\lambda\boxtimes\id$). We get that
$$\ind_{P_I\cap K_n}^{K_n}(\tau\boxtimes\id)\simeq 
\bigcup_{m\geq N_s}\ind_{P_I^0(m)}^{K_n}(\lambda\boxtimes\id).$$
Hence, any typical representation occurs as a sub-representation of
$$\ind_{P_I^0(m)}^{K_n}(\lambda\boxtimes\id),$$ for some  $m\geq N_s$.

We first have to understand the representation
$$\ind^{P_I^0(m)}_{P_I^0(m+1)}(\id),$$
for $m\geq N_s$. It is convenient to define a normal subgroup
$R_I(m)$ of $P^0_I(m)$ such that $P_I^0(m)$ is equal to
$R_I(m)P_I^0(m+1)$, and $R_I(m)\cap P_I^0(m+1)= R_I(m+1)$, for
$m\geq N_s$. For any integer $m\geq N_s$, we define $R_I(m)$ to be the group
$K_{n}(m)\pi^{-1}(K_{n-1}(N_s)\times (1+\mathfrak{P}_F^{N_s}))$. The
group $R_I(m)$ has the Iwahori decomposition with respect to the
parabolic subgroup $P_I$ and its Levi subgroup $M_I$. 
\begin{lemma}
  The group $R_I(m)$ is a normal subgroup of $P_I^0(m)$.  The group $R_I(m+1)$
  is a normal subgroup of $R_I(m)$, for all $m\geq N_s$.
\end{lemma}
\begin{proof}
  By definition of the groups $R_I(m)$, we have
  $R_I(m)\cap U_I=P_I^0(m)\cap U_I$, and
  $R_I(m)\cap \bar{U}_I=P_I^0(m)\cap \bar{U}_I$. 
To show the normality
  of $R_I(m)$ in $P^0_I(m)$, we have to verify that $P_I^0(m)\cap M_I$
  normalize the group $R_I(m)$. But, $P_I^0(m)\cap M_I$ normalizes the
  group $R_I(m)\cap U_I=U_I(\integers{F})$ and
  $R_I(m)\cap \bar{U}_I=\bar{U}_I(\varpi^m_F\integers{F})$. The group
  $K_n(m)\cap M_I$ is a normal subgroup of $M_I(\integers{F})$.
  Hence, $P_I^0(m)\cap M_I$ normalizes $R_I(m)\cap M_I$. This shows the
  first part of the lemma.

  Since $R_I(m)\cap P_I=R_I(m+1)\cap P_I$, we have to check that
  $R_I(m)\cap \bar{U}_I$ normalizes the group $R_I(m+1)$. We note that
  $\bar{U}_I$ is abelian. Hence, we have to check that the conjugations
  $u^{-}j(u^{-})^{-1}$ and $u^{-}u^{+}(u^{-})^{-1}$ belong to the
  group $R_I(m+1)$, for all $u^{-}\in R_I(m)\cap \bar{U}_I$,
  $j\in R_I(m+1)\cap M_I=R_I(m)\cap M_I$, and
  $u^{+}\in R_I(m+1)\cap U_I=U_I(\integers{F})$. Let us begin with the
  element $u^{-}j(u^{-})^{-1}$. We have
  $u^{-}j(u^{-})^{-1}=j\{j^{-1}u^{-}j(u^{-})^{-1}\}$. Let
$$j=\begin{pmatrix}J_1&0\\0&j_1\end{pmatrix}\
 \ \ u^{-}=\begin{pmatrix}1_{n-1}&0\\U^{-}&1\end{pmatrix}$$
be the block diagonal form of $j$ and $u^{-}$;
$J_1\in K_{n-1}(N_s), j_1\in 1+\ideal{F}^{N_s}$ and
$U^{-}\in \varpi_F^m {\rm Mat}_{1\times (n-1)}(\integers{F})$.  The element
$j^{-1}u^{-}j(u^{-})^{-1}$ is of the form
$$\begin{pmatrix}1_{n-1}&&0\\ j_1^{-1}U^{-}J_1-U^{-}&&1 \end{pmatrix}.$$
We note that the matrix $j_1^{-1}U^{-}J_1-U^{-}$ belongs to
$\varpi_F^{m+1}{\rm Mat}_{1\times (n-1)}(\integers{F})$.  This shows
that $j^{-1}u^{-}j(u^{-})^{-1}\in R_I(m+1)\cap \bar{U}_I$. Hence we
get that
$$u^{-}j(u^{-})^{-1}=j\{j^{-1}u^{-}j(u^{-})^{-1}\}\in R_I(m+1).$$

We now consider the conjugation $u^{-}u^{+}(u^{-})^{-1}$. We write
$u^{+}$ in its block matrix form as
$$\begin{pmatrix}1_{n-1}& U^{+}\\0&1\end{pmatrix}$$
where $U^{+}\in {\rm Mat}_{(n-1)\times1}(\integers{F})$. Now the conjugation
$u^{-}u^{+}(u^{-})^{-1}$ in the block matrix from is as follows
$$\begin{pmatrix}1_{n-1}-U^{+}U^{-}&&U^{+}\\ 
  -U^{-}U^{+}U^{-}&&U^{-}U^{+}+1\end{pmatrix}.$$ Since
$U^{-}U^{+}U^{-}\in \varpi^{m+1}{\rm Mat}_{1\times
  (n-1)}(\integers{F})$, we conclude that
$u^{-}u^{+}(u^{-})^{-1}\in R_I(m+1)$. This ends the proof of this
lemma.
\end{proof}
Note that $P_I^0(m)=R_I(m)P_I^0(m+1)$. Using Mackey decomposition, we
get that
$$\res_{R_I(m)}\ind_{P_I^0(m+1)}^{P_I^0(m)}(\id)\simeq \ind^{R_I(m)}_{R_I(m+1)}(\id).$$
It follows from the Iwahori decomposition with respect to $P_I$ and
$M_I$ that the inclusion of $R_I(m)\cap \bar{U}_I$ in $R_I(m)$ induces
an isomorphism between $R_I(m)/R_I(m+1)$ and the abelian group
$$\dfrac{R_I(m)\cap \bar{U}_I}{R_I(m+1)\cap \bar{U}_I}.$$
Hence, the representation $\ind^{R_I(m)}_{R_I(m+1)}(\id)$ decomposes
as a direct sum of characters $\eta_k$, for $1\leq k\leq p$, where
$\eta_k$ is trivial on $R_I(m+1)$.  The group $P_I^0(m)$ acts on these
characters, and let $\{\eta_{n_k}\ |\ n_k\in\{1,2,\dots,p\}\}$ be a
set of representatives for the orbits under this action. Let
$Z(\eta_k)$ be the $P_I^0(m)$-stabiliser of the character $\eta_k$,
for $1\leq k\leq p$.  Now Clifford theory gives us the isomorphism
\begin{equation}\label{equation_n+1_1}
 \ind_{P_I^0(m+1)}^{P_I^0(m)}(\id)=
\bigoplus_{\eta_{n_k}}\ind_{Z(\eta_{n_k})}^{P_I^0(m)}(U_{\eta_{n_k}}),
\end{equation}
where $U_{\eta_{n_k}}$ is any irreducible representation of
$Z(\eta_{n_k})$ such that $\res_{R_I(m)}U_{\eta_{n_k}}$ contains
$\eta_{n_k}$.

We have to bound the group $Z(\eta_k)$. We note that $P_I^0(m)$ is
equal to $(P_I^0(m)\cap M_I)R_I(m)$.  Hence, we have
$Z(\eta_k)=(Z(\eta_k)\cap M_I)R_I(m)$. To bound the group $Z(\eta_k)$ we can
only need to control $Z(\eta_k)\cap M_I$. Let
$u^{-}\in R_I(m)\cap \bar{U}_I$ and
$$\begin{pmatrix}1_{n-1}&0\\U^{-}&1\end{pmatrix}$$
be the block form of $u^{-}$, where $U^{-}$ is a matrix in
$\varpi_F^{m}{\rm Mat}_{1\times (n-1)}(\integers{F})$. The map
$u^{-}\mapsto \varpi_F^{-m}U^{-}$ induces an
$M_I(\integers{F})$-equivariant isomorphism between
${\rm Mat}_{1\times (n-1)}(k_F)$ and the quotient
$$\dfrac{R_I(m)\cap \bar{U}_I}{R_I(m+1)\cap \bar{U}_I}.$$
We also have an $M_I(\integers{F})$-equivariant isomorphism between
${\rm Mat}_{(n-1)\times 1}(k_F)$ and $\widehat{{\rm Mat}_{1\times
    (n-1)}(k_F)}$ (see Lemma \ref{prelim_duality_isom}). We note
that $P_I^0(m)\cap M_I=J^0\times\mathfrak{o}_F^{\times}$.

Let $\eta$ be a non-trivial character of $R_I(m)$ which is trivial on
$R_I(m+1)$. For the present purposes, it is enough to bound the
subgroup $Z(\eta)\cap
(U^0(\mathfrak{B})\times\mathfrak{o}_F^{\times})$, for $\eta\neq
\id$. Since we have a $M_I(\integers{F})$-equivariant isomorphism
between the group of characters on the quotient $R_I(m)/R_I(m+1)$ with
${\rm Mat}_{n-1\times 1}(k_F)$, we can as well study the group
$Z(A)\cap (U^0(\mathfrak{B})\times\mathfrak{o}_F^{\times})$, where
$Z(A)$ is the $M_I(\integers{F})$-stabiliser of a non-zero matrix
$A\in {\rm Mat}_{(n-1)\times 1}(k_F)$. The action of the group
$M_I(\integers{F})$ factorizes through
$K_{n-1}(1)\times(1+\ideal{F})$.  Hence, the group
$(\id_{n-1}+\mathfrak{D}^e)\times (1+\ideal{F})$ is contained in the
kernel of the action of
$U^0(\mathfrak{B})\times\mathfrak{o}_F^{\times}$ obtained by
restriction. Recall that $\mathfrak{D}$ is the radical of
$\mathfrak{B}$.

This reduces our
situation to the following setting. The group
$\g{n-1}{k_F}\times k_F^{\times}$ acts on ${\rm Mat}_{n-1\times 1}(k_F)$ by setting
$$(g_1,g_2)A=g_1Ag_2^{-1},$$
where $g_1\in \g{n-1}{k_F}$, $g_2\in k_F^{\times}$, and
$A\in {\rm Mat}_{n-1\times 1}(k_F)$. The basis
$(v_1, v_2,\dots, v_{n-1})$ we fixed for the vector space $V$ at the
beginning of this section, gives a basis of the $k_F$-vector space
$$(\mathfrak{o}_E/\varpi_F\mathfrak{o}_E)^{(n-1)/ef}=
(\mathfrak{o}_E/\mathfrak{p}_E^e)^{(n-1)/ef}.$$ Such a basis gives the
inclusion
$$\g{(n-1)/ef}{\integers{E}/\ideal{E}^e}=
U^0(\mathfrak{B})/U^e(\mathfrak{B}) \hookrightarrow \g{n-1}{k_F}.$$
Recall that $\lambda=\kappa\otimes\rho$, where $\rho$ is a cuspidal
representation of
$\g{(n-1)/ef}{k_E}=U^0(\mathfrak{B})/U^1(\mathfrak{B})$. Hence, we are
interested in the mod $\ideal{E}$ reduction of the first projection of
$$Z_{\g{(n-1)/ef}{\integers{E}/\ideal{E}^e}\times k_F^{\times}}(A),$$
for some non-zero matrix $A$ in ${\rm Mat}_{(n-1)\times 1}(k_F)$. We
set $n_0=(n-1)/ef$.

Let $\varpi_E$ be an uniformiser of $\integers{E}$. Let $N$ be the
operator on the $k_E$-vector space
$W:=(\mathfrak{o}_E/\mathfrak{p}_E^e)^{n_0}$ given by
$$N(w)=\varpi_E.w\ \text{for all}\ w\in W.$$
Since
$\mathfrak{o}_E/\mathfrak{P}_E^e=k_E\oplus k_E\overline{\varpi_E}\oplus
k_E\overline{\varpi_E}^2\oplus \dots\oplus
k_E\overline{\varpi_E}^{e-1}$,
we obtain a decomposition of $W=W_1\oplus W_2\oplus \dots\oplus W_e$
such that $N$ restricted to $W_i$ is an isomorphism onto $W_{i+1}$ for
$i<e$, and $N$ acts trivially on $W_e$. The mod $\ideal{E}$-reduction
of $W$ is the projection onto the first factor $W_1$.

Any $k_E[N]$-linear map $T$ is determined by its restriction to the
space $W_1$.  Given a map $T\in \ho_{k_E}(W_1,W)$ we obtain an
extension $\tilde{T}\in \End_{k_E[N]}(W)$ by setting
$$\tilde{T}(w)=N^{(i-1)}T(N^{-(i-1)}w),$$ 
for all $w\in W_i$ and $1\leq i\leq e$. The map 
$T\mapsto \tilde{T}$ gives us an isomorphism of vector spaces
\begin{equation}
  \ho_{k_E}(W_1,W)\simeq \End_{k_E[N]}(W,W).
\end{equation}
We may write $W=W_1\oplus NW$. This shows that the mod $\ideal{E}$
reduction map, denoted by $\pi_E$, is given by sending 
$\tilde{T}$ to $p_1\circ\tilde{T}_{|W_1}$, where $p_1$ is the
projection onto the first factor of the direct sum
$W_1\oplus W_2\oplus \dots \oplus W_e$. Now $\End_{k_E}(V_1)$ is a
subspace of $\ho_{k_E}(W_1,W)$ and mod $\ideal{E}$ reduction of
$\widetilde{\End_{k_E}(W_1)}$ (the image of $\End_{k_E}(W_1)$ under
the map $T\mapsto \tilde{T}$) is identity on $\End_{k_E}(W_1)$. Hence
$\Aut_{k_E[N]}(W)$ is the semi-direct product
$\widetilde{\Aut_{k_E}(W_1)}\ker(\pi_E)$.

Let $Q$ be a parabolic subgroup fixing the flag
$\mathcal{F}^i=\oplus_{j=1}^iW_i$, for $1\leq i\leq e$, and $L$ be its
Levi subgroup fixing the decomposition
$W_1\oplus W_2\oplus \dots\oplus W_e$. Now
$\widetilde{\Aut_{k_E}(W_1)}$ diagonally embeds in $L$ and
$\ker(\pi_E)$ is a subgroup of the radical of $Q$. The group
$\g{n-1}{k_F}\times k_F^{\times}$ acts on
${\rm Mat}_{n-1\times 1}(k_F)$ by the map
$(g_1,g_2)A\mapsto g_1Ag_2^{-1}$, where $g_1\in \g{n-1}{k_F}$,
$g_2\in k_F^{\times}$, and $A\in {\rm Mat}_{n-1\times 1}(k_F)$.  We
now have the action of
$\g{n_0}{\mathfrak{o}_E/\mathfrak{P}_E^e}\times k_F^{\times}$ on
${\rm Mat}_{n-1\times 1}(k_F)$ by restriction from the action of
$\g{n-1}{k_F}\times k_F^{\times}$. We are interested in
$$(\pi_E\times\id)\{Z_{\g{n_0}{\mathfrak{o}_E/\mathfrak{P}_E^e}\times k_F^{\times}}(A)\},$$
for some $A\in {\rm Mat}_{n-1\times 1}(k_F)\backslash \{0\}$.  

We first look at $Z_{Q\times k_F^{\times}}(A)$. Let $(A_{ij})$ be an
element of $Q$ in its block form. Let $(A_1,A_2,\dots,A_e)^{\text{T}}$
be the block form of $A$, where $A_j$ is a block of size $1\times
n_0f$. If $k$ is the largest positive integer such that $A_k\neq 0$
and $A_k=0$, then we get that $$A_{kk}A_{k}a^{-1}=A_{k},$$ for all
$((A_{ij}),a)\in Z_{Q\times k_F^{\times}}(A)$.  Hence $\{A_{kk}\ | \
((A_{ij}),a)\in Z_{Q\times k_F^{\times}}(A)\}$ is contained in a
proper parabolic subgroup of $\Aut_{k_F}(W_k)$.  If $n_0f>1$, we get that
$$(\pi_E\times\id)\{Z_{\g{n_0}{\mathfrak{o}_E/\mathfrak{P}_E^e}\times k_F^{\times}}(A)\}$$
is a subgroup of $H\times k_F^{\times}$, where $H$ is a subgroup of
$\Aut_{k_E}(W_1)$ whose image under the inclusion map
$\Aut_{k_E}(W_1)\hookrightarrow \Aut_{k_F}(W_1)$ is contained in a
proper $k_F$-parabolic subgroup of $\Aut_{k_F}(W_1)$, since $A_k\neq
0$.

We recall the following proposition due to Pa\v{s}k\={u}nas (see
\cite[Definition 6.2, lemma 6.5, Proposition
6.8]{Paskunas-uniqueness}).
\begin{proposition}[Pa\v{s}k\={u}nas]\label{proposition_n+1_paskunas}
  Let $W$ be a $k_E$-vector space with (finite) dimension greater than
  one. Let $\rho$ be a cuspidal representation of $\Aut_{k_E}(W)$.
  Let $H$ be a subgroup of $\Aut_{k_E}(W)$ such that the image of $H$
  under the inclusion map $\Aut_{k_E}(W)\hookrightarrow\Aut_{k_F}(W)$
  is contained in a proper parabolic subgroup of $\Aut_{k_F}(W)$. For
  every $H$-irreducible sub-representation $\xi$ of $\res_{H}(\rho)$
  there exists an irreducible representation $\rho'$ of
  $\Aut_{k_E}(W)$ such that $\rho'\not\simeq \rho$ and
  $\ho_{H}(\xi, \rho')\neq 0$.
\end{proposition}

Going back to
$Z(\eta)\cap (U^{0}(\mathfrak{B})\times \mathfrak{o}_F^{\times})$, for
$n_0>1$, we get that for every irreducible sub representation $\xi$ of
$$\res_{Z(\eta)\cap (U^{0}(\mathfrak{B})\times
  \mathfrak{o}_F^{\times})} ((\kappa\otimes\rho)\boxtimes \id),$$ there
exists an irreducible representation $\rho'$ of
$U^0(\mathfrak{B})/U^1(\mathfrak{B})$ such that
$$\ho_{Z(\eta)\cap (U^{0}(\mathfrak{B})\times
  \mathfrak{o}_F^{\times})} (\xi,
(\kappa\otimes\rho')\boxtimes\id)\neq 0.$$

For the case $n_0=1$ and $q_F>2$, we have to look at
\begin{equation}\label{equation_gl_0}
(\pi_E\times\id)\{Z_{{\mathfrak{o}_E/\mathfrak{P}_E^e}^{\times}\times k_F^{\times}}(A)\}
\end{equation}
for some nonzero matrix $A\in {\rm Mat}_{n-1\times 1}(k_F)$.  We notice
that the group (\ref{equation_gl_0}) is of the form
$\{(a,a)|a\in k_F^{\times}\}$ if $k_E=k_F$. Let $k_E$ be a proper
extension of $k_F$. If $(a,b)$ is an element of the centraliser
(\ref{equation_gl_0}) then $aA_kb^{-1}=A_k$ ($A_k$ is defined in the
previous paragraph). This shows that $a$ lies in a proper parabolic
subgroup of $\g{f}{k_F}$. This shows that the group
(\ref{equation_gl_0}) is of the form
$\{(a,b)\ | a\in \mathbb{F}^{\times}, b\in k_F^{\times}\}$ where
$\mathbb{F}$ is a proper sub-field of $k_E$.  

In the case where $n_0=1$ and $k_E=k_F$, we consider a non-trivial
character $\phi$ of
$U^0(\mathfrak{B})/U^1(\mathfrak{B})=k_F^{\times}$. We observe that
$$\res_{Z_{J^0\times
    \mathfrak{o}_F^{\times}}(A)}(\lambda\phi\boxtimes\phi^{-1}) \simeq
\res_{Z_{J^0\times
    \mathfrak{o}_F^{\times}}(A)}(\lambda\boxtimes\id).$$ Moreover,
$[M_I,\sigma\boxtimes\id]$ and $[M_I, \sigma'\boxtimes\phi^{-1}]$ are
two distinct inertial classes for any cuspidal representation
$\sigma'$ containing $(J^0, \lambda\otimes\phi)$. Here,we will use the
same notation $\phi$ for the inflation of $\phi$ to the group
$\mathfrak{o}_F^\times$.

In the case where $n_0=1$ and $k_E$ is a proper extension of $k_F$, we
consider a non-trivial character $\phi$ of $k_E^{\times}$ which is
trivial on $\mathbb{F}^{\times}$. We note that
$$\res_{Z_{J^0\times
    \mathfrak{o}_F^{\times}}(A)}(\lambda\phi\boxtimes\id) \simeq
\res_{Z_{J^0\times \mathfrak{o}_F^{\times}}(A)}(\lambda\boxtimes\id)$$
and moreover $[M_I,\sigma\boxtimes\id]$ and
$[M_I, \sigma'\boxtimes\id]$ are two distinct inertial classes for any
cuspidal representation $\sigma'$ containing
$(J^0, \lambda\otimes\phi)$.  With this we finish our preliminaries.

\subsection{Uniqueness of typical representations}
In this part, we will prove the uniqueness of typical representations
for $[M_I, \sigma\boxtimes\id]$. By Frobenius reciprocity, we get that
$\lambda\boxtimes\id$ occurs with multiplicity one in
$\ind^{P_I^0(N_s)}_{P_I^0(m)}(\lambda\boxtimes\id)$, for all
$m>N_s$.  We denote by $U^0_m(\lambda\boxtimes\id)$ the complement
of $\lambda\boxtimes\id$ in
$\ind^{P_I^0(N_s)}_{P_I^0(m)}(\lambda\boxtimes\id)$. We use the
notation $U_m(\lambda\boxtimes\id)$ for the representation
$$\ind_{P_I^0(N_s)}^{K_n}\{U_m^0(\lambda\boxtimes\id)\}.$$
\begin{theorem}\label{theorem_gl_1}
  Let $\#k_F>2$. The $K_n$-irreducible sub representations of
  $U_m(\lambda\boxtimes\id)$ are atypical, for all $m> N_s$.
\end{theorem}
\begin{proof}
  We prove the theorem by induction on the positive integer $m> N_s$.
  We suppose the theorem is true for some positive integer $m>N_s$. We
  will show the same for $m+1$.

We first note that
$$\ind_{P_I^0(m+1)}^{K_n}(\lambda\boxtimes\id)\simeq 
\ind_{P_I^0(m)}^{K_n}\{\ind_{P_I^0(m+1)}^{P_I^0(m)}(\id)\otimes
(\lambda\boxtimes\id)\}.$$
From the decomposition \ref{equation_n+1_1}, we get that
$$\ind_{P_I^0(m+1)}^{K_n}(\lambda\boxtimes\id)\simeq
\bigoplus_{\eta_{n_k}}\ind_{Z(\eta_{n_k})}^{K_n}
\{(\lambda\boxtimes\id)\otimes U_{\eta_{n_k}}\}.$$ 
Note that the above
sum is taken over the orbits for the action of $P_I(m)$ on the set of
characters of $R_I(m)/R_I(m+1)$. Since there is a unique orbit, among
the characters $\{\eta_k\ |\ 1\leq k\leq p\}$, consisting the identity
character, we get that
\begin{equation}
 \ind_{P_I^0(m+1)}^{K_n}(\lambda\boxtimes\id)\simeq
 \ind_{P_I^0(m)}^{K_n}(\lambda\boxtimes\id)\oplus 
\bigoplus_{\eta_{n_k}\neq
  \id}\ind_{Z(\eta_{n_k})}^{K_n}
\{(\lambda\boxtimes\id)\otimes U_{\eta_{n_k}}\}.
\end{equation}
Let $\Gamma$ be an irreducible sub-representation of 
\begin{equation}
\ind_{Z(\eta_{n_k})}^{K_n}\{(\lambda\boxtimes\id)
\otimes U_{\eta_{n_k}}\}. 
\end{equation}

We have two cases $n_0=1$ and $n_0>1$. If $n_0=1$ we have seen that we
can find a non-trivial character $\phi$ of
$k_E^{\times}=U^0(\mathfrak{B})/U^1(\mathfrak{B})$ such that
$$\ind_{Z(\eta_{n_k})}^{K_n}\{(\lambda\boxtimes\id)
\otimes U_{\eta_{n_k}}\}\simeq
\ind_{Z(\eta_{n_k})}^{K_n}
\{(\lambda\phi\boxtimes\phi^{-1})\otimes U_{\eta_{n_k}}\}$$
or 
$$\ind_{Z(\eta_{n_k})}^{K_n}\{(\lambda\boxtimes\id)\otimes
U_{\eta_{n_k}}\}
\simeq \ind_{Z(\eta_{n_k})}^{K_n}
\{(\lambda\phi\boxtimes\id)\otimes U_{\eta_{n_k}}\}. $$
Hence in this case, the irreducible subrepresentations of
\begin{equation}\label{rev_case_1}
\ind_{Z(\eta_{n_k})}^{K_n}\{(\lambda\boxtimes\id)\otimes U_{\eta_{n_k}}\}
\end{equation}
occur as subrepresentations of
$\res_{K_n}i_{P_I}^{G_n}(\sigma'\boxtimes\chi')$, where $\sigma'$ is a
cuspidal representation of $G_{n-1}$ containing the type $(J^0,
\lambda\otimes\phi))$.  The inertial classes $[G_{n-1}, \sigma]$ and
$[G_{n-1},\sigma']$ are distinct. Hence, any irreducible
subrepresentation of \eqref{rev_case_1} is atypical.

Now consider the case $n_0>1$. In this case, there exists an
irreducible representation $\xi$ of
$(\pi_E\times \id)\{Z(\eta)\cap
(U^0(\mathfrak{B})\times\mathfrak{o}_F^{\times})\}$ such that $\Gamma$
is a sub-representation of
\begin{equation}\label{equation_n+1_2}
\ind_{Z(\eta_{n_k})}^{K_n}\{((\xi\otimes\kappa)\boxtimes\id)\otimes U_{\eta_{n_k}}\}. 
\end{equation}
Now Proposition \ref{proposition_n+1_paskunas} gives us an irreducible
representation $\rho'\not\simeq \rho$ of $U^0(\mathfrak{B})$ obtained
by inflation of an irreducible representation of
$U^0(\mathfrak{B})/U^1(\mathfrak{B})$ such that $\xi$ is contained in
$\rho'$.  Now the representation \eqref{equation_n+1_2} is a
sub-representation of
$$\ind_{Z(\eta_{n_k})}^{K_n}
\{((\rho'\otimes\kappa)\boxtimes\id)\otimes U_{\eta_{n_k}}\}. $$ The
above representation is contained in
\begin{equation}\label{equation_n+1_atypical}
  \ind_{P_I^0(m+1)}^{K_n}((\rho'\otimes\kappa)\boxtimes\id)
  \simeq \ind_{P_I(m+1)}^{K_n}(\tau'\boxtimes\id),
\end{equation}
where $\tau'$ is isomorphic to $\ind_{J^0}^{K_{n-1}}(\rho'\otimes\kappa)$.
The representation $\rho'\otimes\kappa$
is still irreducible (see \cite[Chapter 5, Proposition
5.3.2(3)]{Orrangebook}).

We will show that irreducible subrepresentations of
\eqref{equation_n+1_atypical} are atypical for the inertial class
$[M_I, \sigma\boxtimes\id]$.

Any irreducible sub-representation of (\ref{equation_n+1_atypical})
occurs as a sub-representation of
$$\ind_{P_I(m)}^{K_n}(\gamma\boxtimes\id),$$
where $\gamma$ is a $K_{n-1}$-irreducible subrepresentation of
$\tau'$. Now $\gamma$ is contained in an irreducible smooth
representation say $\sigma_0$ of $G_{n-1}$.  By Frobenius reciprocity this
is possible only if the representation $\rho'\otimes\kappa$ of $J^0$
is contained in $\sigma_0$. We have two possible situations either $\rho'$
is cuspidal or otherwise.  If $\rho'$ is cuspidal, then the
representation $\sigma_0$ is a cuspidal representation such that
$\sigma_0\not\simeq \sigma$. Hence the representation
$$\ind_{P_I(m+1)}^{K_n}(\gamma\boxtimes\id)$$
occurs in 
$$\res_{K_n}i_{P_I}^{G_n}(\sigma_0\boxtimes\id)$$
with $[G_{n-1}\times F^{\times}, \sigma_0\boxtimes\id]\neq
[G_{n-1}\times F^{\times}, \sigma\boxtimes\id].$ This shows that
irreducible subrepresentations of \eqref{equation_n+1_atypical} are
atypical representations.

Consider the case where $\rho'$ is not cuspidal. If
$(J^0, \rho'\otimes\kappa)$ is contained in an smooth irreducible
representation $\sigma_0$, then $\sigma_0$ either contains a
non-maximal simple-type $(J_1^0, \rho_1\otimes\kappa_1)$ or contains a
split type (see \cite[Chapter 8, Theorem 8.3.5]{Orrangebook}). We also refer
to the article \cite[Lemma 2, Proposition 1]{Bush_henniart_intert} for
quick reference. From this we conclude that $\sigma_0$ is not a
cuspidal representation. Hence, the representation
\eqref{equation_n+1_atypical} is contained in
$$\res_{K_n}i_{P'}^{G_n}(\sigma'),$$
where $P'$ is a parabolic subgroup $G_n$ properly contained in $P_I$,
and $\sigma'$ is a cuspidal representation of a Levi subgroup of
$P'$. Since $\mathcal{I}(i_{P'}^{G_n}(\sigma'))\neq [M_I,
\sigma\boxtimes\id]$, we get that the irreducible subrepresentations of
\eqref{equation_n+1_atypical} are atypical.
\end{proof}
Recall the definition of the integer $N_s$ from \eqref{rev_n_1_begin}.
Now, any typical representation for $s$ occurs as a subrepresentation
of $\ind_{P_I^0(m)}^{K_n}(\lambda\boxtimes\id)$, for some $m\geq
N_s$. For $m>N_s$, we have
$$\ind_{P_I^0(m)}^{K_n}(\lambda\boxtimes\id)=
\ind_{P_I^0(N_s)}^{K_n}(\lambda\boxtimes\id)\oplus
U_m(\lambda\boxtimes\id).$$ From the above theorem we get that the
typical representations for $s$ occur as subrepresentations of
$$\ind_{P_I^0(N_s)}^{K_n}(\lambda\boxtimes\id).$$
The above representation may still contain atypical
representations. We will indeed show that this is the case and
complete the classification.

The first observation is that the group $J_s$ in a semisimple
Bushnell-Kutzko type $(J_s, \lambda\boxtimes \id)$, for $s=[M_I,
\sigma\boxtimes\id]$, contains the group $P_I^0(N_s)$. Hence we will
try to decompose the representation
\begin{equation}\label{equation_gl_1}
\ind_{P^0_I(N_s)}^{J_s}(\id).
\end{equation}
We also note that $P_I^0(N_s)\cap P_I=J_s\cap P_I$.  Let $l+1=el'+r$,
where $0\leq r<e$. If $r\leq 1$, then $J_s=P_I^0(N_s)$ and hence, we
have nothing further to analyse and Theorem \ref{theorem_gl_1}
completes the classification of typical representations. {\bf From now
  we assume that $e>2$ and $r>1$}. Note that the depth-zero case is
already handled in \cite{level_zero_gl_n_types} (see Theorem
\ref{thr_depth_zero}).  We will first verify that the group
$U_I(\integers{F})$ acts trivially on the representation
\eqref{equation_gl_1}.

Let $u^{+}$ and $u^{-}$ be two matrices from
$J_s\cap U_I=U_I(\integers{F})$ and $J_s\cap \bar{U}_I$ respectively.
Let $u^{+}$ and $u^{-}$ in block form be written as:
$$\begin{pmatrix}1_{n-1} &U^{+}\\0&1\end{pmatrix}\ \text{and} 
\ \begin{pmatrix}1_{n-1}&0\\U^{-}&1\end{pmatrix}$$ respectively.  The
block form of the conjugation $u^{-}u^{+}(u^{-})^{-1}$ is given by
$$\begin{pmatrix}1_{n-1}-U^{+}U^{-}&&
  U^{+}\\-U^{-}U^{+}U^{-}&& U^{-}U^{+}+1\end{pmatrix}.$$
We have
$$\begin{pmatrix}0&0\\U^{-}&0\end{pmatrix}\in a_{l+1}(\Lambda)\cap
\bar{\mathfrak{n}}_I=\varpi_F^{l'}(a_{r}(\Lambda)\cap
\bar{\mathfrak{n}}_I).$$
 If $r\geq1$, then the valuation of each entry
of a matrix in $a_{r}(\Lambda)\cap \bar{\mathfrak{n}}_I$ is at least
one. This shows that the valuation of each entry in $U^{-}U^{+}U^{-}$
is at least $l'+2$. From which the conjugation
$u^{-}u^{+}(u^{-})^{-1}$ lies in the group $P^0(N_s)$. If $r=0$ and
$l'=0$, then we are in the case where $\sigma$ is a level-zero
cuspidal representation and in this case $J_s=P^0_I(N_s)$. If $r=0$
and $l'>0$, then valuation of each entry in $U^{-}U^{+}U^{-}$ has
valuation $2l'>l'+1$ and hence $u^{-}u^{+}(u^{-})^{-1}\in P_I^0(N_s)$.
Hence, the group $U_I(\integers{F})$ acts trivially on the
representation \eqref{equation_gl_1}.

From the Iwahori decomposition of the group $J_s$, we get that $J_s$
is equal to $(J_s\cap \bar{P}_I)P_I^0(N_s)$. Hence we have:
$$\res_{J_s\cap \bar{P}_I}\ind_{P_I^0(N_s)}^{J_s}(\id)\simeq 
\ind^{J_s\cap \bar{P}_I}_{P^0_I(N_s)\cap \bar{P}_I}(\id).$$ Note that
$J_s\cap \bar{P}_I$ is a semi-direct product of the groups
$(J_s\cap M_I)$ and $(J_s\cap \bar{U}_I)$. Let $\eta_k$, for
$1\leq k\leq t$, be all the characters of the group
$J_s\cap \bar{U}_I$ which are trivial on the group
$P_I^0(N_s)\cap \bar{U}_I$. The group $J_s\cap \bar{P}_I$ acts on
these characters. Let $\{\eta_{k_p}\}$ be a set of representatives for
the orbits under this action. We denote by $Z(\eta_{k_p})$ the
$J_s\cap \bar{P}_I$ stabiliser of the character $\eta_{k_p}$. Let
 $U_{\eta_{k_p}}$ be the isotypic component of the
character $\eta_{k_p}$ in the representation
$$\ind^{J_s\cap \bar{P}_I}_{P^0_I(N_s)\cap \bar{P}_I}(\id).$$
The space $U_{\eta_{k_p}}$ has a natural action of
$Z(\eta_{k_p})$. Now Clifford theory gives the decomposition
$$ \ind^{J_s\cap \bar{P}_I}_{P^0_I(N_s)\cap \bar{P}_I}(\id)\simeq
\bigoplus_{\eta_{k_p}} 
\ind^{J_s\cap \bar{P}_I}_{Z(\eta_{k_p})}(U_{\eta_{k_p}}).$$
We note that the character $\id$ occurs with a
multiplicity one in the list of characters $\eta_k$.

If $K_s$ is the kernel of the representation \eqref{equation_gl_1},
then $K_s\cap Z(\eta_{k_p})$ acts trivially on $U_{\eta_{k_p}}$. Hence
we can extend the representation $U_{\eta_{k_p}}$ to the group
$Z(\eta_{k_P})K_s$ such that $K_s$ acts trivially on the
extension. Now consider the representation
$$\pi=\ind^{J_s}_{Z(\eta_{k_P})K_s}U_{\eta_{k_p}}.$$
Note that $K_s\cap \bar{P}_I$ is contained in the group
$Z(\eta_{k_p})\cap \bar{P}_I$ and moreover, $U_I(\integers{F})$ is
contained in $K_s$. Hence we have $J_s=(J_s\cap \bar{P}_I)Z(\eta_{k_p})K_s$.
From Mackey decomposition, we have
$$\res_{J_s\cap \bar{P}_I}
\ind^{J_s}_{Z(\eta_{k_p})K_s}U_{\eta_{k_p}}\simeq 
\ind^{J_s\cap \bar{P}_I}_{Z(\eta_{k_p})K_s\cap (J_s\cap
  \bar{P}_I)}(U_{\eta_{k_p}})
\simeq \ind^{J_s\cap \bar{P}_I}_{Z(\eta_{k_p})}(U_{\eta_{k_p}}).$$
We hence have
\begin{equation}\label{equation_gl_2}
  \ind^{J_s}_{P_I^0(N_s)}(\id)\simeq 
  \bigoplus_{\eta_{k_p}}\ind^{J_s}_{Z(\eta_{k_p})K_s}U_{\eta_{k_p}}.
\end{equation}

Now using the decomposition (\ref{equation_gl_2}) we get  the
decomposition
$$\ind_{P^0_I(N_s)}^{K_n}(\lambda\boxtimes\id)
\simeq \bigoplus_{\eta_{k_p}}\ind^{K_n}_{Z(\eta_{k_p})K_s}
\{U_{\eta_{k_p}}\otimes(\lambda\boxtimes\id)\}.$$ Note that the
character $\id$ occurs with multiplicity one among the characters
$\eta_k$. Moreover, we have $Z(\id)K_s=(J_s\cap
\bar{P}_I)K_s=J_s$ and we get that
\begin{equation}
  \ind_{P^0_I(N_s)}^{K_n}(\lambda\boxtimes\id)
  \simeq  
  \ind_{J_s}^{K_n}(\lambda\boxtimes\id)\oplus 
  \bigoplus_{\eta_{k_p}\neq
    \id}\ind^{K_n}_{Z(\eta_{k_p})K_s}
  \{U_{\eta_{k_p}}\otimes(\lambda\boxtimes\id)\}.
\end{equation}
\begin{lemma}\label{lemma_gl_2}
  Let $\#k_F>2$ and $\eta_{k_p}$ be a non-trivial character. The
  irreducible sub representations of
$$\ind^{K_n}_{Z(\eta_{k_p})K_s}
\{U_{\eta_{k_p}}\otimes(\lambda\boxtimes\id)\}$$
are atypical.
\end{lemma}
\begin{proof}
  We observe that
  $Z(\eta_{k_p})=(Z(\eta_{k_p})\cap M_I)(J_s\cap \bar{U}_I)$. This
  shows that we have to bound the group $Z(\eta_{k_P})\cap M_I$, for
  $\eta_{k_p}\neq \id$. Recall that $\eta_k$, for $1\leq k\leq t$, are
  the characters of the quotient group
\begin{equation}\label{equation_gl_3}
\dfrac{(J_s\cap \bar{U}_I)}{(P_I^0(N_s)\cap \bar{U}_I)}.
\end{equation} 
Now let $u^{-}$ be a matrix from the group $J_s\cap \bar{U}_I$. In the
block form the matrix $u^{-}$ is of the form
$$\begin{pmatrix}1_{n-1}&0\\U^{-}&1\end{pmatrix}$$
where $U^{-}=[M_1, M_2, \dots, M_{e}]$, $M_i$ is a matrix of size
$(1\times (n-1)/e)$. Let $\delta=N_s-1$. The map $\Phi$
$$[M_1, M_2, \dots, M_{e}]\mapsto [\varpi_F^{\delta}M_1, 
\varpi_F^{\delta}M_2, \dots, \varpi_F^{\delta}M_e]$$ identifies the
quotient \eqref{equation_gl_3} with a subspace $\mathfrak{t}_1$ of
${\rm Mat}_{1\times n-1}(k_F)$. This identification commutes with the
action of $M_I\cap J_s$, since $\Phi$ is none other than conjugation
by an element from the $Z(M_I)$ (The centre of $M_I$).

Let $\mathfrak{t}_2$ be the following space of column matrices:
$$\mathfrak{t}_2=\{(0, 0, \dots, 0, M_{r},\dots, M_{e})^{\text{T}}\ |
M_j\in {\rm Mat}_{(n-1)/e\times 1}(k_F)\ \forall\ r\leq j\leq e\}.$$
The group $M_I(\mathfrak{o}_F)\simeq K_{n-1}\times K_1$ acts on the
space ${\rm Mat}_{(n-1)\times 1}(k_F)$ via the conjugation action of $M_I$ on
$U_I$. The space $\mathfrak{t}_2$ is stable under the action of
$U^0(\mathfrak{A})\times\mathfrak{o}_F^\times\subseteq M_I$. The 
pairing $XY$, where $X\in\mathfrak{t}_2$ and
$Y\in \mathfrak{t}_1$, gives a perfect pairing between
$\mathfrak{t}_1$ and $\mathfrak{t}_2$. This pairing is equivariant for
the action of
$U^0(\mathfrak{A})\times\mathfrak{o}_F^\times\subseteq M_I$. This
gives an identification of the space of characters of $\mathfrak{t}_1$
with the space $\mathfrak{t}_2$ in a
$U^0(\mathfrak{A})\times\mathfrak{o}_F^\times$ equivariant way.

The group $M_I(\mathfrak{o}_F)$ acts on
${\rm Mat}_{(n-1)\times 1}(k_F)$ through its quotient $M_I(k_F)$. Now,
the action of the group
$(U^0(\mathfrak{B})\times \mathfrak{o}_F^{\times})\subseteq J_s\cap
M_I$ on $\mathfrak{t}_2$ factors through its quotient by its subgroup
$(1_n+\mathfrak{D}^e)\times (1+\ideal{F})$. Let $A$ be a non-zero
matrix in $\mathfrak{t}_2$.

Recall that $n_0=(n-1)/ef$. 
Now recall that we denote by $\pi_E$ by mod $\ideal{E}$ reduction map.
We have seen that (the paragraph above the proposition
\ref{proposition_n+1_paskunas})
$$(\pi_E\times\id)\{Z_{\g{n_0}{\mathfrak{o}_E/\mathfrak{P}_E^e}\times k_F^{\times}}(A)\}$$
is a subgroup of $H\times k_F^{\times}$, where $H$ is a subgroup of
$\g{n_0}{k_E}$ whose image under the inclusion map
$\g{n_0}{k_E}\hookrightarrow \g{n-1}{k_F}$ is contained in a proper
$k_F$-parabolic subgroup of $\g{n-1}{k_F}$. From the result of
Pa\v{s}k\={u}nas, stated as Proposition
\ref{proposition_n+1_paskunas}, we get that for every irreducible
representation $\xi$ of
$$\res_{Z(\eta_{k_p})}\{U_{k_p}\otimes((\kappa\otimes\rho)\boxtimes\id)\}$$
we can find an irreducible representation $\rho'\not\simeq \rho$ such
that $\xi$ occurs in the representation
$$\res_{Z(\eta_{k_p})}\{U_{k_p}\otimes((\kappa\otimes\rho')\boxtimes\id)\}.$$
Hence irreducible subrepresentations of
$$\ind^{K_n}_{Z(\eta_{k_p})K_s}
\{U_{\eta_{k_p}}\otimes(\lambda\boxtimes\id)\}$$
occur as  subrepresentations of 
$$\ind^{K_n}_{Z(\eta_{k_p})K_s}
\{U_{\eta_{k_p}}\otimes((\kappa\otimes\rho')\boxtimes\id)\}.$$
Now the above representation occurs as a sub-representation of 
$$\ind^{K_n}_{P^0_I(N_s)}\{(\kappa\otimes\rho')\boxtimes\id\}
\simeq \ind^{K_n}_{P_I(N_s)}(\tau'\boxtimes\id)\},$$ where $\tau'$ is
given by
$$\ind_{J^0}^{K_{n-1}}(\kappa\otimes\rho').$$

Any irreducible sub representation $\gamma$ of $\tau'$ occurs in an
irreducible smooth representation $\sigma_0$ of $\g{n-1}{F}$. Assume that
$\rho'$ is cuspidal. The representation $\kappa\otimes\rho'$ is
contained in the representation $\res_{J^0}\gamma$ and hence is
contained in $\sigma_0$. This implies that $\sigma_0$ is cuspidal but not
inertially equivalent to $\sigma$. If $\rho'$ is not cuspidal, then
the representation $\sigma_0$ is not cuspidal. Hence, in every case, $\sigma_0$
is not inertially equivalent to $\sigma$. This shows that irreducible
sub representations of
$$ \ind^{K_n}_{P_I(N_s)}(\gamma\boxtimes\id)$$
are atypical. We conclude the lemma. 
\end{proof}
\begin{theorem}\label{rev_n+1_final}
  Let $n>2$ and $q_F>2$. Let $\Gamma$ be any typical representation
  for the inertial class $s=[M_I, \sigma\boxtimes\chi]$. The representation
  $\Gamma$ is isomorphic to the representation
$$\ind_{J_s}^{K_n}(\lambda_s),$$
where $(J_s, \lambda_s)$ is any Bushnell-Kutzko semi-simple type for
the inertial class $s$. If $P$ is a parabolic subgroup containing $M_I$ as
a Levi subgroup, then $\Gamma$ occurs with a multiplicity one in the
representation
$$\res_{K_n}i_P^{G_n}(\sigma\boxtimes\chi).$$
\end{theorem}
\begin{proof}
  Let $G$ be the group of $F$-rational points of a connected reductive
  group defined over $F$.  For any inertial class $t=[L, \Theta]$ of
  $G$, and for $g\in N_G(L)$ we define $t^g$ to be the inertial class
  $[L, \Theta^g]$. The map sending $t$ to $t^g$ is well defined. We
  denote by $N_G(t)$ the group $\{g\in N_G(L)\ |\ t^g=t\}$.  The group
  $N_G(t)$ clearly contains the group $L$ and the quotient
  $W_t=N_G(t)/L$ is finite. The cardinality of $W_t$ does not depend
  on the choice of $L$. We return to the case where $G=G_n$ and
  $t=s$. The intertwining of the representation
  $\ind_{J_s}^{K_n}(\lambda_s)$ is bounded by the cardinality of
  $W_s$. We have $|W_s|=1$ since $n>2$. Hence, the representation
  $\ind_{J_s}^{K_n}(\lambda_s)$ is irreducible. We refer to
  \cite{Smoothrepcptgr}[Lemma 11.5] for these results. Hence the
  uniqueness of the typical representation. The multiplicity follows
  from the results \ref{theorem_gl_1} and \ref{lemma_gl_2}.
\end{proof}
\section{Principal series components}
Let $I$ be the partition $(1,1,\dots,1)$ of $n$. Recall that we denote
by $B_n$ the group $P_I$, $U_n$ the group $U_I$, and $T_n$ the group
$M_I$ respectively. In this section, we will classify typical
representations for the inertial classes $s=[T_n, \chi]$, where $\chi$
is a character of $T_n$. Let $\tau$ be a typical representation for
the inertial class $s$.  The compact induction $\ind_{K_n}^{G_n}\tau$
is a finitely generated representation. Let $\pi$ be an irreducible
quotient of $\ind_{K_n}^{G_n}\tau$. By Frobenius reciprocity, the
$K_n$-representation $\tau$ occurs in the $G_n$ representation $\pi$.

Let $B$ be any Borel subgroup of $G_n$, $T$ be a maximal split torus
of $G_n$ contained in $B$, and $\chi'$ be a character of $T$. If
$(T, \chi')$ and $(T_n, \chi)$ are inertially equivalent, then the
representation $\pi$ occurs as a sub-quotient of
$i_{B}^{G_n}(\chi'')$, where $\chi''$ is obtained from $\chi'$ by
twisting with an unramified character of $T$. For classifying typical
representations it is enough to say which $K_n$-irreducible sub
representations of $i_{B}^{G_n}(\chi')$ are typical for the inertial class
$[T_n, \chi]$.

Let $\sigma$ be a permutation of the set $\{1,2,\dots,n\}$. Let
$\chi=\boxtimes_{i=1}^{n}\chi_i$ be any character of
$T_n=\prod_{i=1}^n{F^{\times}}$.  We denote by $\chi^{\sigma}$ the
character $\boxtimes_{i=1}^{n}\chi_{\sigma(i)}$ of $T_n$.  We observe
that the pairs $(T_n, \chi^\sigma)$ and $(T_n, \chi)$ are inertially
equivalent. This implies that for a classifying typical representations
we can classify typical representation occurring in
$i_{B_n}^{G_n}\chi^\sigma$, for any $\sigma$. We
will use a convenient permutation $\sigma$ which satisfies the
condition in the following lemma.
\begin{lemma}\label{lemma_princ_1}
  Given any sequence of characters $x_i=\chi_i$ of
  $\mathfrak{o}_F^{\times}$, there exists a permutation
  $\{y_i \ |1\leq i\leq n\}$ of $\{x_i\ |1\leq i\leq n\}$ such that
$$l(y_iy_k^{-1})\geq \max\{l(y_iy_j^{-1}), l(y_jy_k^{-1})\},$$
for all $1\leq i\leq j\leq k\leq n$.
\end{lemma}
\begin{proof}
  For any ultrametric space $(X,d)$ and given any $n$ points
  $x_1,x_2,x_3,\dots,x_n$ in $X$ we may choose a permutation
  $y_1,y_2,\dots,y_n$ of the sequence $\{x_i|1\leq i\leq n\}$ such
  that
$$d(y_i,y_k)\geq \max\{d(y_i,y_j),d(y_j,y_k)\}$$
for all $i\leq j\leq k$. Now apply this fact to the space $X$
consisting of characters of $\mathfrak{o}_F^{\times}$ and the distance
function $d(\chi_1,\chi_2)$ is defined as the level
$l(\chi_1\chi_2^{-1})$ if $\chi_1\neq \chi_2$ and $0$ otherwise.  We
point out that this ordering is not unique in general. We refer to
\cite[Lemma 1] {howe_principal_series}for a proof of these results.
\end{proof}
\begin{remark}
 We note that the condition 
 $l(y_iy_k^{-1})\geq \max\{l(y_iy_j^{-1}), l(y_jy_k^{-1})\}$
 is equivalent to an equality since we always have 
 $$l(y_iy_k^{-1})\leq \max\{l(y_iy_j^{-1}), l(y_jy_k^{-1})\}.$$
\end{remark}
Given an inertial class $[T_n, \chi]$ we choose the
representative $(T_n, \chi^{\sigma})$ where $\sigma$ is a permutation
such that
$$l(\chi_{\sigma(i)}\chi_{\sigma(k)}^{-1})\geq 
\max\{l(\chi_{\sigma(i)}\chi_{\sigma(j)}^{-1}),
l(\chi_{\sigma(j)}\chi_{\sigma(k)}^{-1})\}.$$ 
{\bf From now on we
  assume that the pair $(T_n, \boxtimes_{i=1}^{n}\chi_i)$ satisfies
  the condition }
\begin{equation}\label{equation_princ_ineq}
l(\chi_{i}\chi_{k}^{-1})\geq \max\{l(\chi_{i}\chi_{j}^{-1}),
l(\chi_{j}\chi_{k}^{-1})\},
\end{equation}
for all $i\leq j\leq k$.

In the following subsection we construct subgroups $H_m$, for $m\geq
1$ such that 
\begin{enumerate}
\item $H_1=J_s$, where $J_s$ is the compact open subgroup of a
  Bushnell--Kutzko type $(J_s, \chi)$ of $s$,
\item $H_{m+1}\subset H_{m}$, for all $m\geq 1$ and
  $\bigcap_{m\geq 1}H_m=B_n\cap K_n$,
\item The representation $\chi$ of $T_n\cap K_n$ extends to a
  representation of $H_m$ such that $H_m\cap U_n$ and $H_m\cap
  \bar{U}_n$ are contained in the kernel of this extension.
\end{enumerate}
Such a construction gives the following equality:
$$\ind_{K_n\cap B_n}^{K_n}\chi=\bigcup_{m\geq
  1}\ind_{H_m}^{K_n}\chi.$$ Later, we show that any $K_n$-irreducible
sub representation of $\ind_{H_{m+1}}^{K_n}\chi /\ind_{H_m}^{K_n}\chi$
is atypical.
\subsection{Construction of compact open 
subgroups \texorpdfstring{$H_m$}{}}
Let $\mathcal{A}=(a_{ij})$ be a lower nilpotent matrix of size
$n\times n$ such that $a_{ij}$ is non-negative, for $i> j$, and
\begin{equation}\label{imp_princ}
a_{ki}=\max\{a_{ji}, a_{kj}\},
\end{equation} 
for $1\leq i< j< k\leq n$.  We denote
by $J(\mathcal{A})$ the set of $n\times n$ matrices $(m_{pq})$ such
that $m_{pq}\in \integers{F}$, for $p<q$, and
$m_{pq}\in \mathfrak{P}_F^{a_{pq}}$, for $p\geq q$.  As a consequence
of the condition $a_{ki}=\max\{a_{ji}, a_{kj}\}$ we get two important
inequalities
\begin{equation}\label{key_ineq_1}
a_{i1}\geq a_{i2}\geq\dots\geq a_{ii-1}
\end{equation}
and 
\begin{equation}\label{key_ineq_2}
a_{j+1j}\leq a_{j+2j}\leq \dots\leq a_{nj}.
\end{equation}
The first is a consequence of $a_{ik-1}=\max\{a_{kk-1}, a_{ik}\}$, for
$k<i$, and the second is a consequence of
$a_{k+1j}=\max\{a_{k+1k}, a_{kj}\}$, for $j<k$.
\begin{lemma}\label{lemma_princ_2}
The set $\mathcal{J}(\mathcal{A})$ is an order in 
${\rm Mat}_{n\times n}(\integers{F})$ 
\end{lemma}
\begin{proof}
  The set $\mathcal{J}(\mathcal{A})$ is an additive group. We now
  check that the set $\mathcal{J}(\mathcal{A})$ is closed under
  multiplication. Let $(m_{ij})$ and $(m_{ij}')$ be two matrices from
  $J(\mathcal{A})$. If $i>j$, then the $i\times j$ term in the product
  matrix $(m_{ij})(m_{ij}')$ can be split into three terms:
  $$t_1:=m_{i1}m_{1j}'+m_{i2}m_{2j}'+\dots+m_{ij}m_{ji}',$$
  $$t_2:=m_{ij+1}m_{j+1k}'+\dots+m_{ii}m_{ij}'$$ and
  $$t_3:=m_{ii+1}m_{i+1j}'+\dots+m_{in}m'_{nj}.$$ Observe that
  $\nu_F(m_{ik}m_{kj}') \geq a_{i1}$, for $k\leq j$.  This shows
  that $\nu_F(t_1)\geq \min\{a_{i1}, a_{i2},\dots ,a_{ij}\}$ and
  $$\min\{a_{i1},\dots ,a_{ij}\}\geq a_{ij}.$$ The valuation 
  $\nu_F(m_{ik}m_{kj}')\geq a_{ik}+a_{kj}$, for all $j\leq k\leq i$,
  and $a_{ik}+a_{kj}$ is greater or equal to $a_{ij}$.  We get that
  $\nu_F(t_2)\geq a_{ij}$. Finally the valuation
  $\nu_F(m_{ik}m_{kj})\geq a_{kj}$, for $k>i$. The valuation
  $\nu_F(t_3)\geq \min\{a_{i+1j},\dots, a_{nj}\}$ and
  $\min\{a_{i+1j},\dots, a_{nj}\}\geq a_{ij}$. Hence the additive
  group $\mathcal{J}(\mathcal{A})$ is closed under
  multiplication. Since $\mathcal{J}(\mathcal{A})$ is an
  $\integers{F}$ lattice in ${\rm Mat}_{n\times n}(F)$ we get that
  $\mathcal{J}(\mathcal{A})$ is an order in
  ${\rm Mat}_{n\times n}(\integers{F})$.
\end{proof}

We denote by $J(\mathcal{A})$ the set of invertible elements of
$\mathcal{J}(\mathcal{A})$. The following are examples of
$J(\mathcal{A})$.
\begin{enumerate}
\item If $\mathcal{A}=0$ then the group $J(\mathcal{A})$ is
  $K_n$.
\item If $\mathcal{A}=(a_{ij})$ with $a_{ij}=1$, for $i>j$, then
  $J(\mathcal{A})$ is the Iwahori subgroup with respect to the
  standard Borel subgroup $B_n$.
\end{enumerate}
The examples $(1)$ and $(2)$ satisfy Iwahori decomposition with
respect to the standard Borel subgroup $B_n$. The next lemma concerns
the Iwahori decomposition of $J(\mathcal{A})$ in general.

Let $\mathcal{A}=(a_{ij})$ be a lower nilpotent matrix such that
$a_{ki}=\max\{ a_{kj},a_{ji}\}$, for $1\leq i< j< k\leq n$. We define
an ordered partition $I$ of $n$ by induction on the set of positive
integers $m\leq n$.  Let $I_1:=(1)$ and if we know
$I_m=(n_1,n_2,\dots,n_r)$, for some $m\leq n-1$, then $I_{m+1}$ is the
following partition 
\begin{eqnarray}
  I_{m+1}=\begin{cases}
    (n_1,n_2,\dots,n_r,1)\ \text{if}\ a_{m+1m}\neq 0\\
    (n_1,n_2,\dots,n_r+1)\ \text{otherwise}.
\end{cases}
\end{eqnarray}
We denote by $I(\mathcal{A})$ the partition $I_n$.

\begin{lemma}\label{lemma_princ_3}
  The group $J(\mathcal{A})$ satisfies Iwahori decomposition with
  respect to the parabolic subgroup $P_{I(\mathcal{A})}$ and the Levi
  subgroup $M_{I(\mathcal{A})}$. We have
  $J(\mathcal{A})\cap
  M_{I(\mathcal{A})}=M_{I(\mathcal{A})}(\integers{F})$,
  $J(\mathcal{A})\cap
  U_{I(\mathcal{A})}=U_{I(\mathcal{A})}(\integers{F})$.
\end{lemma}

\begin{proof}
  We use induction on the positive integer $n$. If $n=1$ then
  $J(\mathcal{A})$ is $\mathfrak{o}_F^{\times}$ and the lemma is
  vacuously true. We assume that the lemma is true for all positive
  integers less than $n$. Let $I(\mathcal{A})$ be the ordered
  partition $(n_1,n_{2},\dots n_r)$. If $r=1$ then the lemma is true
  by default. We suppose $r>1$. We will show below that every element
  $j\in J(\mathcal{A})$ can be written as a product $u_1j_1$ with
  $u_1\in J(\mathcal{A})\cap \bar{U}_{(n_1,n-n_1)}$ and
  $j_1\in J(\mathcal{A})\cap P_{(n_1,n-n_1)}$.  Now $j_1$ can be
  written as $j_2u_1^{+}$ where
  $u_1^{+}\in U_{(n_1,n-n_1)}(\integers{F})$ and
  $j_2\in M_{(n_1,n-n_1)}\cap J(\mathcal{A})$.  Now $j_2$ can be
  written as $j_3u_2^{+}$ where
  $j_3\in J(\mathcal{A}) \cap M_{I(\mathcal{A})}$ and
  $u_2^{+}\in U_{I(\mathcal{A})}(\integers{F})$.  The group
  $J(\mathcal{A})\cap M_{(n_1,n-n_1)}$ is equal to $K_{n_1}\times
  J(\mathcal{A}')$
  where the nilpotent matrix $\mathcal{A}'=(a_{ij}')$ is given by
  $a_{ij}'=a_{i+n_1j+n_1}$. By induction hypothesis $J(\mathcal{A}')$
  satisfies Iwahori decomposition with respect to the standard
  parabolic subgroup $P_{I(\mathcal{A}')}$ and its Levi subgroup
  $M_{I(\mathcal{A}')}$ and
  $I(\mathcal{A}')=(n_{2},n_{3},\dots,n_r)$. Let $j_3=(j_3^0, j_3^1)$
  where $j_3^0\in K_{n_1}$ and
  $j_3^1\in J(\mathcal{A}')$. Now $j_3^1=u_3^{-}j_4u_3^{+}$ where
  $u_3^{-}\in \bar{U}_{I(\mathcal{A}')}\cap J(\mathcal{A}')$,
  $u_3^{+}\in {U}_{I(\mathcal{A}')}\cap J(\mathcal{A}')$ and
  $j_4\in M_{I(\mathcal{A}')}\cap J(\mathcal{A}')$. Hence
  $j=u_1u_3^{-}(j_3^0,j_4)u_3^{+}u_2^{+}$ (with a slight abuse of
  notation the elements $u_3^{-}$ and $u_3^{+}$ are considered as
  elements of $\bar{U}_{I(\mathcal{A})}$ and $U_{I(\mathcal{A})}$
  respectively and $(j_3^0,j_4)$ is an element of
  $J(\mathcal{A})\cap M_{I(\mathcal{A})}=K_{n_1}\times
  (J(\mathcal{A}')\cap M_{I(\mathcal{A}')})$).

  We now prove that $j\in J(\mathcal{A})$ can be written as a product
  $u_1j_1$ with $u_1\in J(\mathcal{A})\cap \bar{U}_{(n_1,n-n_1)}$ and
  $j_1\in J(\mathcal{A})\cap P_{(n_1,n-n_1)}$. Let $j=(j_{pq})$.  Let
  $C_i^1$ be the $i^{th}$-column of the first diagonal block (of size
  $n_1\times n_1$) on the diagonal. If every entry of $C_i^1$ has
  positive valuation then, we claim that the all the entries of the
  $i^{th}$ column $C_i$ have positive valuation.  Suppose the $k^{th}$
  entry $j_{ki}$ of $C_i$ is an unit for some $k>n_1$. This shows that
  $a_{ki}$ the $ki^{th}$-entry of $\mathcal{A}$ is zero. Now the
  inequality (\ref{key_ineq_1}) gives $a_{ki}\geq a_{kn_1}$ and this
  implies that $a_{kn_1}=0$. Now note that $a_{kn_1}\geq a_{n_1+1n_1}$
  from the inequality (\ref{key_ineq_2}).  This shows that
  $a_{n_1+1n_1}$ is zero which gives a contradiction from the
  definition of $I(\mathcal{A})$.  We now deduce that $j_{ki}$ is not
  invertible. This shows the claim. Since $j$ is invertible we
  conclude that at least one entry of $C_i^1$ is an unit. Let
  $E_{ij}(c)=I_n+e_{ij}(c)$ where $e_{ij}(c)$ is the matrix with its
  $ij$ entry $c$ and all other entries $0$.  The left multiplication
  of $E_{ij}(c)$ results in the row operation $R_j+cR_i$. Since at
  least one entry of $C_i^1$ is an unit we assume that its
  $q^{th}$-entry is an unit. We can perform row operations $R_p+cR_q$
  for all $p\geq n_1$ to make the $p^{th}$-entry trivial. We also note
  that the elementary matrix corresponding to this row operation also
  belongs to the group $J(\mathcal{A})$ (note that $q\leq n_1\leq
  p$). This completes the task of making $j$ as the product $u_1j_1$.
The uniqueness of the Iwahori decomposition is standard. 
\end{proof} 

Let $s=[T_n, \chi]$ be an inertial equivalence class. Let $m$ be a
positive integer and $\mathcal{A}_{\chi}(m)$ be the lower nilpotent
matrix $(a_{ij}^m)\in {\rm Mat}_{n\times n}(\mathbb{Z})$, where
$$a_{ij}^m=l(\chi_i\chi_j^{-1})+m-1,$$
for $n\geq i>j\geq 1$.  As
shown earlier, the representative
$(T_n, \chi=\boxtimes_{i=1}^n\chi_i)$ for $s$, can be chosen such
that
$$a_{ik}=\max\{a_{ij}, a_{jk}\},$$ 
for all $i<j<k$. We denote by $J_{\chi}(m)$ the group
$J(\mathcal{A}_{\chi}(m))$. Note that
$J_{\chi}(m')\subseteq J_{\chi}(m)$, for all $m'\geq m$.  In our
situation we have $I(\mathcal{A}_{\chi}(m))$ is $(1,1,\dots,1)$, since
none of $a_{ii+1}^m$ are zero. Hence, using lemma \ref{lemma_princ_3}
the group $J_{\chi}(m)$ satisfies the Iwahori decomposition with
respect to $B_n$ and $T_n$.

\begin{lemma}\label{lemma_princ_4_1}
  The character $\chi=\boxtimes_{i=1}^n\chi_i$ of $T_n\cap K_n$
  extends to a character of $J_{\chi}(1)$ such that
  $J_{\chi}(1)\cap U_n$ and $J_{\chi}(1)\cap \bar{U}_n$ are contained
  in the kernel of the extension.
\end{lemma}
\begin{proof}
  Let $m=(m_{ij})$ be an element of $J_{\chi}(1)$. We define
  $\tilde{\chi}(m)=\prod_{i=1}^n\chi_i(m_{ii})$. The verification that
  $\tilde{\chi}$ is a character of the group $J_{\chi}(1)$ is very
  computational in nature and we sketch the proof here and for
  complete details see \cite[Section 3, Lemma 3.1, Lemma 3.2
  ]{AllanRoche} or \cite[Pg 278-279]{howe_principal_series}.  The idea
  is to get an open normal subgroup $U$ of $J_{\chi}(1)$ such that
  $J_{\chi}(1)/U$ is isomorphic to $T(\integers{F})/T_{\chi}$ where
  $T_{\chi}$ is an open subgroup of $T(\integers{F})$ which is
  contained in the kernel of $\chi$. The subgroup $U$ is generated by
  $J_{\chi}(1)\cap \bar{U}_n$ and
  $J_{\chi}(1) \cap U_n=U_n(\integers{F})$. One shows that $U$
  satisfies Iwahori decomposition with respect to the Borel subgroup
  $B_n$ and $U\cap T_n$ is given by
  $\prod_{\alpha\in
    \Phi}\alpha^{\vee}(1+\mathfrak{P}_F^{l(\chi\alpha^{\vee})})$ where
  $\Phi$ is the set of roots of $\gl_n$ with respect to $T_n$ and
  $\alpha^{\vee}$ stands for the dual root. We observe that
  $U\cap T_n$ is contained in the kernel of $\chi$.
 \end{proof}
 We are now ready to define the sequence of groups $H_m$. We set
 $H_m=J_{\chi}(m)$, for $m\geq 1$. We now get the equality
$$\res_{K_n}i_{B_n}^{G_n}(\chi)=\bigcup_{m\geq 1}\ind_{J_{\chi}(m)}^{K_n}(\chi).$$
For the purposes of proofs by induction, we need to construct some
more compact open subgroups of $K_n$.

We denote by $\mathcal{A}_{\chi}(1,m)$ the lower nilpotent matrix
$(a_{ij})$ where $a_{ij}=l(\chi_i\chi_j^{-1})$ for $j<i<n$,
$a_{nj}=l(\chi_n\chi_j^{-1})+m-1$, for $1\leq j<n-1$. Given a lower
nilpotent matrix $\mathcal{A}=(a_{ij})$ such that
$a_{ki}=\max\{a_{kj}, a_{ji}\}$ we associated a compact subgroup
$J(\mathcal{A})$. The matrix $\mathcal{A}_{\chi}(1,m)$ need not
satisfy this condition but, we can still associate the group
$J(\mathcal{A}_{\chi}(1,m))$ to the matrix $\mathcal{A}_{\chi}(1,m)$.
We will prove this in the next Lemma.

\begin{lemma}\label{lemma_princ_4}
  Let $\mathcal{J}(\mathcal{A}_{\chi}(1,m))$ be the set consisting of
  matrices $(m_{ij})\in {\rm Mat}_{n\times n}(\mathfrak{o}_F)$ such
  that $m_{ij}\in \mathfrak{P}_F^{a_{ij}}$ for all $i,j$.  The set
  $\mathcal{J}(\mathcal{A}_{\chi}(1,m))$ is an order in
  ${\rm Mat}_{n\times n}(\integers{F})$.
\end{lemma}
\begin{proof}
  The set $\mathcal{J}(\mathcal{A}_{\chi}(1,m))$ is a lattice
  in ${\rm Mat}_{n\times n}(F)$ and we have to verify that
  $\mathcal{J}(\mathcal{A}_{\chi}(1,m))$ is closed under
  multiplication.  Let $(m_{ij})$ and $(m_{ij}')$ be two elements of
  the set $\mathcal{J}(\mathcal{A}_{\chi}(1,m))$. We suppose
  $i>j$. The $ij^{th}$-term of the product $(m_{ij})(m_{ij}')$ is the
  sum of the terms:
 $$t_1:=m_{i1}m_{1j}'+m_{i2}m_{2j}'+\dots+m_{ij}m_{ji}',$$
 $$t_2:=m_{ij+1}m_{j+1k}'+\dots+m_{ii}m_{ij}'$$
  and 
  $$t_3:=m_{ii+1}m_{i+1j}'+\dots+m_{in}m'_{nj}.$$
  Note that
  $$\nu_F(t_1)\geq \min\{\nu_F(m_{ik}m_{kj}')\ |\ \text{for all}\ 1\leq k\leq
  j\},$$ and $\nu_F(m_{ik}m_{kj}')=a_{ik}$. If
  $i<n$, then $a_{ik}=l(\chi_{i}\chi_{k}^{-1})$ and
  $a_{ij}\leq a_{ik}$, for all $k\leq j<i$.  This shows that
  $\nu_F(t_1)\geq a_{ij}$. If $i=n$, then we have
  $$a_{ik}=a_{nk}=l(\chi_n\chi_{k}^{-1})+m-1\geq
  l(\chi_i\chi_j^{-1})+m-1=a_{ij},$$ for all $k\leq j<n$. We conclude
  that in every possibility $\nu_F(t_1)\geq a_{ij}$.
  
Consider the term $t_2$. We have 
$$\nu_F(t_2)\geq
\{\nu_F(m_{ik}m_{kj}')\ |\ \text{for all}\ j<k\leq i\}$$ and
$\nu_F(m_{ik}m_{kj}')=a_{ik}+a_{kj}$, for $j<k\leq i$. If $i< n$, then
we have $a_{ik}=l(\chi_i\chi_k^{-1})$ and
$a_{kj}=l(\chi_k\chi_j^{-1})$. From our assumption on the arrangement
of characters $\chi_i$, for $1\leq i\leq n$, we get that
  $$l(\chi_i\chi_j^{-1})=\max\{l(\chi_i\chi_k^{-1}),
  l(\chi_k\chi_j^{-1})\}.$$
  At the same time $i<n$ implies that $a_{ij}=l(\chi_i\chi_j^{-1})$. This
  shows that $\nu_F(m_{ik}m_{kj}')$ is equal to
  $a_{ik}+a_{kj}$ and $a_{ik}+a_{kj} \geq a_{ij}$.  Consider the case
  $i=n$ and in this case,
  $a_{nk}=l(\chi_n\chi_k^{-1})+m-1$. Now $a_{kj}=l(\chi_k\chi_j^{-1})$
  and $a_{nj}=l(\chi_n\chi_j^{-1})+m-1$. From the equality
  $l(\chi_i\chi_j^{-1})=\max\{l(\chi_i\chi_k^{-1}),
  l(\chi_k\chi_j^{-1})\}$ we deduce that
$$l(\chi_i\chi_k^{-1})+l(\chi_k\chi_j^{-1})>l(\chi_i\chi_j^{-1}).$$
Now, adding $m-1$ on both sides we get $a_{ij}>a_{ik}+a_{kj}$. We
conclude that  $\nu_F(t_2)\geq a_{ij}$.

We observe that $\nu_F(m_{ik}m_{kj}')$, for $i<k<n$ is equal to
$a_{kj}=l(\chi_{k}\chi_j^{-1})$ and we have
$l(\chi_{k}\chi_j^{-1}) \geq a_{ij}$. Note that
$a_{nj}=l(\chi_{n}\chi_j^{-1})+m-1\geq a_{ij}$ from which we conclude
that $\nu_F(t_3)\geq a_{ij}$. This shows that
$\nu_F(t_1+t_2+t_3)\geq a_{ij}$ and we prove our lemma.
\end{proof}

Let $J_{\chi}(1,m)$ be the group of units of
$\mathcal{J}(\mathcal{A}_{\chi}(1,m))$. We will need the structure of
the representation
$$\ind_{J_{\chi}(1,m+1)}^{J_{\chi}(1,m)}(\id).$$
We follow a similar strategy to that from previous section.  Let $(a_{ij})$ be
the matrix $\mathcal{A}_{\chi}(1,m)$. Let $K_{\chi}(1,m)$ be the set
of matrices $(m_{ij})$ such that $m_{ij}\in \ideal{F}$, for $i<j<n$,
$m_{in}\in \integers{F}$, for $i<n$, $m_{ii}\in 1+\ideal{F}$, for
$i\leq n$, $m_{ij}\in \ideal{F}^{a_{ij}}$, for $i>j$. In the block
form the group $K_\chi(1,m)$ is given by:
$$\begin{pmatrix}\mathcal{A}_\chi(1)\cap K_{n-1}(1)&{\rm Mat}_{(n-1)\times
    1}(\mathfrak{o}_F)\\
  \mathfrak{n}&1+\mathfrak{p}_F\end{pmatrix},$$ where $\mathfrak{n}$
is the lattice
$(\mathfrak{p}_F^{l(\chi_1\chi_n^{-1})}, \dots,
\mathfrak{p}_F^{l(\chi_{n-1}\chi_n^{-1})})$.
\begin{lemma}\label{lemma_princ_5}
 The set $K_{\chi}(1,m)$ is a normal subgroup of $J_{\chi}(1,m)$.
\end{lemma}
\begin{proof}
  We first check that $K_{\chi}(1,m)$ is closed under matrix
  multiplication. Let $(m_{ij})$ and $(m_{ij}')$ be two matrices from
  the set $K_{\chi}(1,m)$.  Let $i<j<n$ the $ij^{th}$ term is the sum
  of
$$t_1=m_{i1}m_{1j}'+m_{i2}m_{2j}'+\dots+m_{ii}m_{ij}',$$
 $$t_2=m_{ii+1}m_{i+1j}'+m_{ii+2}m_{i+2j}'+\dots+m_{ij}m_{jj}'$$
  and 
 $$t_3=m_{ij+1}m_{j+1j}'+\dots+m_{in}m_{nj}'.$$
  
 Observe that $\nu_F(m_{kj}')>0$, for $1<k\leq i$ and hence,
 $\nu_F(t_1)>0$. Now $\nu_F(m_{ik})>0$, for $i<k\leq j$ and we get
 that $\nu_F(t_2)>0$. Note that $\nu_F(m_{kj}')>0$, for $j<k\leq n$
 and hence, the valuation $\nu_F(t_3)>0$.  This shows that
 $ij^{th}$-term of the matrix product has positive valuation. The
 verifications on congruence conditions for the $in^{\rm th}$-term are
 exactly the same as in Lemma \ref{lemma_princ_4}. The existence of
 inverse for an element in $K_{\chi}(1,m)$ follows from Gaussian
 elimination.

 Now we establish the normality of $K_{\chi}(1,m)$. The group
 $K_{\chi}(1,m)$ satisfies the Iwahori decomposition with respect to the
 subgroups $P_{(n-1,1)}$ and $M_{(n-1,1)}$.  We also note that
 $K_{\chi}(1,m)\cap U_{(n-1,1)}$ is equal to
 $J_{\chi}(1,m)\cap U_{(n-1,1)}$ and
 $K_{\chi}(1,m)\cap \bar{U}_{(n-1,1)}$ is equal to
 $J_{\chi}(1,m)\cap \bar{U}_{(n-1,1)}$. To check the normality of
 $K_{\chi}(1,m)$ we have to check that $J_{\chi}(1,m)\cap M_{(n-1,1)}$
 normalizes $K_{\chi}(1,m)$. This is equivalent to checking that
 $K_{\chi}(1,m) \cap M_{(n-1,1)}$ is a normal subgroup of
 $J_{\chi}(1,m)\cap M_{(n-1,1)}$.
 
We note that
 $J_{\chi}(1,m)\cap M_{(n-1,1)}=J_{\chi'}(1)\times
 \mathfrak{o}_F^{\times}$
 where $\chi'=\boxtimes_{i=1}^{n-1}\chi_i$. Let $p_1$ be the
 projection of $J_{\chi}(1,m)\cap M_{(n-1,1)}$ onto $J_{\chi'}(1)$ and
 $\pi_1$ be the reduction mod $\ideal{F}$ map. Note that
 $K_{\chi}(1,m)\cap M_{(n-1,1)}$ is the kernel of $\pi_1\circ p_1$.
\end{proof}

From the above lemma, the group $K_{\chi}(1,m)$ is  a normal subgroup
of $J_{\chi}(1,m)$. We also note that 
$J_{\chi}(1,m)\cap \bar{U}_{(n-1,1)}$ is contained in $K_{\chi}(1,m)$.
From this we conclude that
$J_{\chi}(1,m)=K_{\chi}(1,m)J_{\chi}(1,m+1)$. From the Mackey
decomposition we get that
$$\res_{K_{\chi}(1,m)}
\ind_{J_{\chi}(1,m+1)}^{J_{\chi}(1,m)}(\id)\simeq 
\ind_{K_{\chi}(1,m)\cap J_{\chi}(1,m+1)}^{K_{\chi}(1,m)}(\id).$$
From the definition of $K_{\chi}(1,m)$ we get that 
$K_{\chi}(1,m)\cap J_{\chi}(1,m+1)=K_{\chi}(1,m+1)$ and 
\begin{equation}\label{equation_princ_1} 
  \res_{K_{\chi}(1,m)}\ind_{J_{\chi}(1,m+1)}^{J_{\chi}(1,m)}(\id)
  \simeq \ind_{K_{\chi}(1,m+1)}^{K_{\chi}(1,m)}(\id).
\end{equation}
\begin{lemma}\label{lemma_princ_6}
 The group $K_{\chi}(1,m+1)$ is a normal subgroup of $K_{\chi}(1,m)$.
\end{lemma}

\begin{proof}
  The group $K_{\chi}(1,m)$ has the Iwahori decomposition with respect
  to $P_{(n-1, 1)}$ and $M_{(n-1, 1)}$,
  $K_{\chi}(1,m)\cap U_{(n-1,1)}$ is equal to
  $K_{\chi}(1,m+1)\cap U_{(n-1,1)}$ and
  $K_{\chi}(1,m)\cap M_{(n-1,1)}$ is equal to
  $K_{\chi}(1,m+1)\cap M_{(n-1,1)}$. We have to check that
  $u^{-}j(u^{-})^{-1}$ and $u^{-}u^{+}(u^{-})^{-1}$ belong to
  $K_{\chi}(1,m+1)$, for
  $u^{-} \in K_{\chi}(1,m)\cap \bar{U}_{(n-1,1)}$,
  $j \in K_{\chi}(1,m)\cap M_{(n-1,1)}$ and
  $u^{+} \in K_{\chi}(1,m)\cap U_{(n-1,1)}$.

  We first consider the case $u^{-}j(u^{-})^{-1}$. We can rewrite
  $u^{-}j(u^{-})^{-1}$ as $j\{j^{-1}u^{-}j(u^{-})^{-1}\}$. Since
  $j\in K_{\chi}(1,m)\cap M_{(n-1,1)}=K_{\chi}(1,m+1)\cap
  M_{(n-1,1)}$, it is enough to show that $j^{-1}u^{-}j(u^{-})^{-1}$
  belongs to the group $K_{\chi}(1,m+1)$. Let $j$ and $u^{-}$ be
  written in their block matrix form as follows.
$$j=
\begin{pmatrix}
J_1&0\\
0&j_1 \end{pmatrix} 
\\ 
\ \ 
u^{-}=
\begin{pmatrix}
1_{n-1}&0\\
U^{-}&1
\end{pmatrix}
$$
The conjugation $j^{-1}u^{-}j(u^{-})^{-1}$ in its block form is given
by
$$\begin{pmatrix}
1_{n-1} &&0\\
j_1^{-1}U^{-}J_1-U^{-}&&1
\end{pmatrix}$$
Let $U^{-}=[u_1,u_2,\dots,u_{n-1}]$ and $J_1=(j_{ij})$. The $k^{th}$
entry of the matrix $U^{-}J_1$ is the sum of
$$t_1=u_1j_{1k}+u_{2}j_{2k}+\dots+u_{k-1}j_{k-1k},$$
$$t_2=u_kj_{kk},$$
 and
$$t_3=u_{k+1}j_{k+1k}+\dots +u_{n-1}j_{n-1k}.$$  If
$l(\chi_k\chi_n^{-1})>1$, then valuation $\nu_F(u_tj_{tk})$, for $t<k$, is
at least $l(\chi_t\chi_n^{-1})+m-1+1\geq
l(\chi_k\chi_n^{-1})+m$. Hence, we have 
$$\nu_F(t_1)\geq l(\chi_k\chi_n^{-1})+m-1.$$ The valuation $\nu_F(u_ta_{tk})$, for
$k<t$, is at least
$l(\chi_t\chi_n^{-1})+l(\chi_t\chi_{k}^{-1})+m-1>
l(\chi_{k}\chi_{n}^{-1})+m-1$.  This shows that
$t_1+t_2+t_3\equiv t_2=u_kj_{kk}=u_k \mod
\ideal{F}^{l(\chi_k\chi_n^{-1})+m}$.  We note that
$j_1^{-1}u-u\in\mathfrak{p}_F^{x+1}$, for any $u\in \mathfrak{p}_F^x$.
hence the matrix
$$\begin{pmatrix}1_{n-1}&&0\\ j_1^{-1}U^{-}J_1-U^{-}&& 1 \end{pmatrix} $$
is contained in $K_{\chi}(1,m+1)\cap \bar{U}_{(n-1,1)}$

Let us consider the conjugation $u^{-}u^{+}(u^{-})^{-1}$. We write
$u^{+}$ in the block form as
$$\begin{pmatrix}
1_{n-1}&U^{+}\\ 
0& 1
\end{pmatrix}$$
The conjugated matrix $u^{-}u^{+}(u^{-})^{-1}$ is given by 
$$\begin{pmatrix}
  1_{n-1}-U^{+}U^{-} && U^{+} \\
  -U^{-}U^{+}U^{-} && U^{-}U^{+}+1\end{pmatrix}.$$ Let
$1_{n-1}-U^{+}U^{-}=(u_{ij})$.  The valuation
$\nu_F(u_{ij})\geq l(\chi_n\chi_j^{-1})$, for $i>j$, and
$l(\chi_n\chi_j^{-1})$ is greater or equal to
$l(\chi_i\chi_j^{-1})$. From this we conclude that
$u^{-}u^{+}(u^{-})^{-1}\in K_{\chi}(1,m+1)$.
\end{proof}

\subsection{Calculation of some stabilisers}
The inclusion map of $K_{\chi}(1,m)\cap \bar{U}_n$ in $K_{\chi}(1,m)$
induces an isomorphism of the quotient $K_{\chi}(1,m)/K_{\chi}(1,m+1)$
with the abelian group
\begin{equation}\label{equation_princ_2}
 \dfrac{K_{\chi}(1,m)\cap \bar{U}_{(n-1,1)}}
{K_{\chi}(1,m+1)\cap \bar{U}_{(n-1,1)}}.
\end{equation}
Hence the representation $\ind_{K_{\chi}(1,m+1)}^{K_{\chi}(1,m)}(\id)$
splits into direct sum of characters say
$\{\eta_k\ |\ \text{for}\ 1\leq k\leq p\}$. The group $J_{\chi}(1,m)$
acts on these characters and let $Z(\eta_k)$ be the
$J_{\chi}(1,m)$-stabiliser of the character $\eta_k$.  From Clifford
theory we get that
\begin{equation}\label{equation_princ_3}
\ind_{J_{\chi}(1,m+1)}^{J_{\chi}(1,m)}(\id)\simeq 
\bigoplus_{\eta_{n_k}}\ind_{Z(\eta_{n_k})}^{J_{\chi}(1,m)}(U_{\eta_{n_k}}),
\end{equation}
where $\eta_{n_k}$ is a representative for an orbit under the action
of $J_{\chi}(1,m)$ and $U_{\chi_{n_k}}$ is an irreducible
representation of the group $Z(\eta_{n_k})$.  Since
$$J_{\chi}(1,m)=(J_{\chi}(1,m)\cap M_{(n-1,1)})K_{\chi}(1,m)$$ 
we get that $Z(\eta_k)=(Z(\eta_k)\cap M_{(n-1,1)})K_{\chi}(1,m)$.

The final step in our preliminaries is to understand the mod
$\ideal{F}$ reduction of the group
$$Z(\eta_k)\cap M_{(n-1,1)}$$
for some non-trivial character $\eta_k$.  The group
$J_{\chi}(1,m)\cap M_{(n-1,1)}$ (which is
$J_{\chi'}(1)\times\mathfrak{o}_F^{\times}$, for
$\chi'=\boxtimes_{i=1}^{n-1}\chi_i$) acts on the quotient
 \begin{equation}\label{equation_princ_33333}
 \dfrac{K_{\chi}(1,m)\cap \bar{U}_{(n-1,1)}}{K_{\chi}(1,m+1)\cap \bar{U}_{(n-1,1)}}
 \end{equation}
 by conjugation. Let $j$ and $u^{-}$ be two elements from
 $J_{\chi}(1,m)\cap M_{(n-1,1)}$ and
 $K_{\chi}(1,m)\cap \bar{U}$ respectively. We write the elements $j$ and
 $u^{-}$ written in their block diagonal form as
 $$\begin{pmatrix}J_1&0\\0&j_1\end{pmatrix} 
\ \ \text{and} \ \ \begin{pmatrix}1_{n-1}&0\\U^{-}&1\end{pmatrix}$$
 respectively. The map $u^{-}\mapsto \varpi_F^{-(m-1)}U^{-}$ induces
 an isomorphism between the group \eqref{equation_princ_33333} and
 ${\rm Mat}_{1\times (n-1)}(k_F)$. 

 The map $u^{-}\mapsto \varpi_F^{-(m-1)}U^{-}$ gives an
 $J_{\chi}(1)\times \integers{F}^{\times}$-equivariant map between
 $M_{1\times (n-1)}(k_F)$ and the group
 (\ref{equation_princ_33333}). We also have a
 $M_{(n-1,1)}(k_F)$-equivariant map between the group of characters of
 $M_{1\times (n-1)}(k_F)$ and the group $M_{(n-1)\times 1}(k_F)$ (see
 Lemma \ref{prelim_duality_isom}).  Hence we obtain a
 $J_{\chi}(1)\times \integers{F}^{\times}$ equivariant map between the
 group of characters of the quotient (\ref{equation_princ_33333}) and
 the group $M_{(n-1)\times 1}(k_F)$, where $J_{\chi}(1)$ acts through
 its (mod $\ideal{F}$) quotient $B_{(n-1)}(k_F)\times k_F^{\times}$
 and the action is $(b, x)A=bAx^{-1}$. Hence to understand the group
 $Z(\eta_k)\cap M_{(n-1,1)}$ for non-trivial $\eta_k$, we first look
 at $Z_{B_{n-1}(k_F)\times k_F^{\times}}(A)$ for some non-zero matrix
 $A$ in $M_{(n-1)\times 1}(k_F)$ .
 
Let $p$ be the projection of $B_{n-1}(k_F)\times k_F^{\times}$ onto
 the diagonal torus
$$T_{n-1}(k_F)\times k_F^{\times}=T_{n}(k_F),$$ let $p_i$ be the
$i^{th}$ projection of $T_n(k_F)$ onto $k_F^{\times}$. The centraliser
$Z_{B_{n-1}(k_F)\times k_F^{\times}}(A)$ of a non-zero matrix
$A=[u_1,u_2,\dots,u_{n-1}]^{\text{T}}$ satisfies the following
property: there exists a $j<n$ such that $p_j(p(t))=p_n(p(t))$, for
all $t\in Z_{B_{n-1}(k_F)\times k_F^{\times}}(A)$ (see \cite[Lemma
3.8]{level_zero_gl_n_types}).  This shows that for any non-trivial
character $\eta_{n_k}$, $Z(\eta_{n_k})\cap T_n$ satisfies the property
that
$$p_j(t)\equiv p_n(t)\ \text{mod}\  \ideal{F}.$$

The character $\chi=\boxtimes_{i=1}^n\chi_i$ of $J_{\chi}(1)$ occurs
with  multiplicity one in the representation
$$\ind_{J_{\chi}(m)}^{J_{\chi}(1)}(\chi).$$
We denote by $U_m^0(\chi)$ the complement of $\chi$ in
$\ind_{J_{\chi}(m)}^{J_{\chi}(1)}(\chi)$. We denote by $U_m(\chi)$ the
representation
$$\ind_{J_{\chi}(1)}^{K_n}\{U_m^0(\chi)\}.$$

\subsection{Elimination of atypical representations}
\begin{theorem}
  Let $q_F>3$. The irreducible sub representations of $U_m(\chi)$
  are atypical. If $n=3$ and $k_F>2$, then the irreducible sub
  representations of $U_m(\chi)$ are atypical.
\end{theorem}

\begin{proof}
  We prove the theorem by using induction on the positive integers $n$
  and $m$.  For $n=1$ the representation $U_m(\chi)$ is trivial and
  the theorem is vacuously true.  Let $n$ be a positive integer
  greater than one. We assume that the theorem is proved for all
  positive integers less than $n$.  We will use the induction
  hypothesis to show the theorem for $n$. 

  We note that $J_{\chi}(1,m)$ and $J_{\chi}(m)$ satisfy the Iwahori
  decomposition with respect to the parabolic subgroup $P_{(n-1,1)}$
  and its Levi subgroup $M_{(n-1,1)}$; we have
  $J_{\chi}(1,m)\cap U_{(n-1,1)}=J_{\chi}(m)\cap U_{(n-1,1)}$ and
  $J_{\chi}(1,m)\cap \bar{U}_{(n-1,1)}=J_{\chi}(m)\cap
  \bar{U}_{(n-1,1)}$.  Hence, the representation
  $\ind^{J_{\chi}(1,m)\cap M_{(n-1, 1)}}_{J_{\chi}(m)\cap
    M_{(n-1,1)}}(\chi)$ extends to a representation of $J_{\chi}(1,m)$
  and this extension is given by
$$\ind_{J_{\chi}(m)}^{J_{\chi}(1,m)}(\chi).$$
If we denote by $\chi'$ the character $T_{n-1}=\boxtimes_{i=1}^{n-1}\chi_i$ of
$\prod_{i=1}^{n-1}F^{\times}$, then we have
$$\ind^{J_{\chi}(1,m)\cap M_{(n-1,1)}}_{J_{\chi}(m)\cap M_{(n-1,1)}}(\chi)
\simeq \ind^{J_{\chi'}(1)}_{J_{\chi'}(m)}(\chi')\boxtimes\chi_n.$$
We also have 
$$\ind^{J_{\chi'}(1)}_{J_{\chi'}(m)}(\chi')\boxtimes \chi_n 
\simeq U_m^0(\chi')\boxtimes\chi_n\oplus \chi.$$
Combining the above isomorphisms we get that 
\begin{equation}
  \ind_{J_{\chi}(m)}^{K_n}(\chi) \simeq
  \ind_{J_{\chi}(1,m)}^{K_n}\{U_m^0(\chi')\boxtimes\chi_n\}
  \bigoplus \ind_{J_{\chi}(1,m)}^{K_n}(\chi).
 \end{equation}
 
 We will use the  induction hypothesis to show that
 $K_n$-irreducible sub representations of
\begin{equation}\label{equation_princ1112}
\ind_{J_{\chi}(1,m)}^{K_n}\{U_m^0(\chi')\boxtimes\chi_n\} 
\end{equation}
are atypical representations. By induction hypothesis any
$K_{n-1}$-irreducible sub-representation of
$U_m(\chi')$ occurs as sub-representation of some
$$i_{P_I}^{G_{n-1}}(\sigma)$$
where $[T_{n-1}, \chi']$ and $[M_I, \sigma]$ are two distinct inertial
classes. We now get that irreducible sub representations of
\eqref{equation_princ1112} occur as sub representations of
$$i_{P_{I'}}^{G_n}(\sigma\boxtimes\chi_n)$$
where $I'$ is obtained from $I$ by adding $1$ at the end of the
ordered partition $I$ of $n-1$. If $I\neq (1,1,\dots,1)$ then the Levi
sub-groups $M_{I'}$ and $T_n$ are not conjugate and hence the inertial
classes $[M_{I'}, \sigma\boxtimes\chi_n]$ and $[T_n, \chi]$ are
distinct inertial classes and this proves our claim in this case. 

Now, we assume that $M_I=T_{n-1}$ and
$\sigma=\boxtimes_{i=1}^{n-1}\sigma_i$ be the tensor factorisation of
the character $\sigma$ of $T_{n-1}$. Since the inertial classes
$[T_{n-1}, \chi']$ and $[T_{n-1}, \sigma]$ are distinct we get a
character $\chi_t$ occurring with non-zero multiplicity in the
multi-set $\{\chi_1,\chi_2,\dots,\chi_{n-1}\}$ but with a different
multiplicity in the multi-set
$\{\sigma_1, \sigma_2,\dots,\sigma_{n-1}\}$.  Adding the character
$\chi_n$ to both multi-sets above keeps the multiplicities of the
character $\chi_t$ distinct and this shows that $[T_n, \chi]$ and
$[T_n, \sigma\boxtimes\chi_n]$ are different inertial classes.

This shows that any typical representation must occur as a
sub-representation of
$$\ind_{J_{\chi}(1,m)}^{K_n}(\chi).$$
The character $\chi$ occurs with  multiplicity one in the
representation $\ind_{J_{\chi}(1,m)}^{J_{\chi}(1)}(\chi)$. We denote
by $U_{1,m}^0(\chi)$ the complement of the character $\chi$ in
$\ind_{J_{\chi}(1,m)}^{J_{\chi}(1)}(\chi)$. We denote by
$U_{1,m}(\chi)$ the representation
$$\ind_{J_{\chi}(1)}^{K_n}\{U_{1,m}^0(\chi)\}.$$
We first note that 
\begin{equation}\label{princ_rev_1}
U_m(\chi)\simeq
\ind_{J_{\chi}(1,m)}^{K_n}\{U_m^0(\chi')\boxtimes\chi_n\} 
\oplus U_{1,m}(\chi).
\end{equation}
We already showed that the $K_n$-irreducible sub representations of
the first summand on the right-hand side of the equation
\eqref{princ_rev_1} are atypical. We now show that $K_n$-irreducible sub
representations of $U_{1,m}(\chi)$ are atypical and this proves the
main theorem.

We first note that 
$$\ind_{J_{\chi}(1,m+1)}^{J_{\chi}(1)}(\chi)\simeq
 \ind_{J_{\chi}(1,m)}^{J_{\chi}(1)}
\{\ind_{J_{\chi}(1,m+1)}^{J_{\chi}(1,m)}(\id)\otimes\chi\}.$$
Using the decomposition \eqref{equation_princ_3} we get that 
$$\ind_{J_{\chi}(1,m+1)}^{J_{\chi}(1)}(\chi)\simeq 
\bigoplus_{\eta_{n_k}}\ind^{J_{\chi}(1)}_{Z(\eta_{n_k})}\{U_{\eta_{n_k}}\otimes\chi\}.$$
Recall that $\eta_{n_k}$ is a representative for the orbit under the
action of the group $J_{\chi}(1,m)$ on the set of characters
$\{\eta_k|\ 1\leq k\leq p\}$ of the group
$K_{\chi}(1,m)/K_{\chi}(1,m+1)$, and $Z(\eta_{n_k})$ is the
$J_{\chi}(1,m)$-stabiliser of the character $\eta_{n_k}$. There is
exactly one orbit consisting of the identity character and hence
\begin{equation}
  \ind_{J_{\chi}(1,m+1)}^{J_{\chi}(1)}(\chi)\simeq
  \ind_{J_{\chi}(1,m)}^{J_{\chi}(1)}(\chi)
  \bigoplus_{\eta_{n_k}\neq \id}
  \ind^{J_{\chi}(1)}_{Z(\eta_{n_k})}\{U_{\eta_{n_k}}\otimes\chi\}. 
\end{equation}
 
Consider the representation
$$\ind^{J_{\chi}(1)}_{Z(\eta_{n_k})}\{U_{\eta_{n_k}}\otimes\chi\}$$
for some representative $\eta_{n_k}\neq \id$. Now, recall that
$Z(\eta_{n_k})\cap T_n$ is a subgroup of
$T_n(\integers{F})=\prod_{i=1}^n\integers{F}^{\times}$, and there
exists a positive integer $j<n$ such that $p_j(t)\equiv p_n(t)$
mod $\ideal{F}$, for all $t\in Z(\eta_{n_k})$. Let $\kappa$ be a
character of $F^{\times}$ such that $\kappa$ is ramified and
$1+\mathfrak{P}_F$ is contained in the kernel of $\kappa$. Let
$\chi^{\kappa}$ be the character
$$\chi_1\boxtimes\chi_2\boxtimes\chi_j\kappa
\boxtimes\dots\boxtimes\chi_n\kappa^{-1}.$$
We observe that $\res_{Z(\eta_{n_k})}(\chi)=
\res_{Z(\eta_{n_k})}(\chi^{\kappa})$ and hence
\begin{equation}\label{equation_princ_compare}
\ind^{J_{\chi}(1)}_{Z(\eta_{n_k})}\{U_{\eta_{n_k}}\otimes\chi\}\simeq 
\ind^{J_{\chi}(1)}_{Z(\eta_{n_k})}\{U_{\eta_{n_k}}\otimes\chi^{\kappa}\}. 
\end{equation}

From the above paragraph we get that
$$U_{1,m+1}^0(\chi)\simeq U_{1,m}^0(\chi)\bigoplus_{\eta_{n_k}\neq
  \id} 
\ind^{J_{\chi}(1)}_{Z(\eta_{n_k})}\{U_{\eta_{n_k}}\otimes\chi\}.$$
and from the above identity we conclude that
\begin{equation}\label{equation_princ_res}
U_{1,m+1}(\chi)\simeq U_{1,m}(\chi)\bigoplus_{\eta_{n_k}\neq \id} 
\ind^{K_n}_{Z(\eta_{n_k})}\{U_{\eta_{n_k}}\otimes\chi\}. 
\end{equation}

From the equation (\ref{equation_princ_compare}) we get that
$$\ind^{K_n}_{Z(\eta_{n_k})}\{U_{\eta_{n_k}}\otimes\chi\}\simeq 
\ind^{K_n}_{Z(\eta_{n_k})}\{U_{\eta_{n_k}}\otimes\chi^{\kappa}\}.$$
If we choose $\kappa$ such that $[T_n, \chi]$ and
$[T_n, \chi^{\kappa}]$ are two distinct inertial classes, then we can
conclude that irreducible sub representations of
$$\ind^{K_n}_{Z(\eta_{n_k})}\{U_{\eta_{n_k}}\otimes\chi\}$$
are atypical. Hence, using the identity (\ref{equation_princ_res})
recursively we get that irreducible sub representations of
$U_{1,m}(\chi)$ are atypical representations, for all positive
integers $m$.

To prove the theorem we have to justify that we can choose a character
$\kappa$ as in the previous paragraph.  Now for any character $\kappa$
non-trivial on $\integers{F}^{\times}$ (such a character exists since
$q_F>2$) and trivial on $1+\ideal{F}$, the equality of the inertial
classes $[T_n, \chi]$ and $[T_n, \chi^{\kappa}]$ implies the
equality of multiplicities of $\chi_j$ in the multi-sets
$\{\chi_1,\chi_2,\dots,\chi_n\}$ and
$\{\chi_1, \chi_2,\dots, \chi_j\kappa,\dots,\chi_n\kappa^{-1}\}$. The
equality of multiplicities implies $\chi_j\chi_n^{-1}=\kappa$.  If
$q_F>3$, then we have at least two non-trivial tame characters and hence
we can choose $\kappa$ distinct from a possibly tame character
$\chi_j\chi_n^{-1}$.

If $l(\chi_i)=1$, for $1\leq i\leq n-1$, and $l(\chi_n)>1$, then we
can always find $\kappa$ such that $[T, \chi]$ and $[T, \chi^\kappa]$
are distinct inertial classes. We note that the induction hypothesis
here is supplied by depth-zero case stated as Theorem \ref{thr_depth_zero}.

Consider the case where $n=3$, $q_F=3$ and $\eta$ is the non-trivial
character of $k_F^\times$. We have the character
$\chi=\chi_1\boxtimes\chi_2\boxtimes \chi_3$ of $T_3$. Assume that
there exists $i\neq j$ and $i,j\in\{1,2,3\}$ such that
$\chi_i\chi_j^{-1}=\eta$, with $l(\eta)=1$. If such a pair $(i, j)$
does not exist, then our present proof goes through. Now, twisting
with the character $\chi_j$ if necessary, and permuting the characters
$\chi_1$, $\chi_2$ and $\chi_3$ if necessary, we may assume that
$\chi_1=\id$, $\chi_2=\eta$. This arrangement still satisfies the
condition \eqref{equation_princ_ineq}. If $l(\chi_3)=1$, then we are
depth-zero case and we refer to Theorem \ref{thr_depth_zero} for a
proof of this result. If $l(\chi_3)>1$, then $K_2$ irreducible
subrepresentations of $U_m(\chi_1\boxtimes\chi_2)$ are atypical, for
$m\geq 1$ (we refer to \ref{thr_depth_zero} or
\cite[Appendix]{Henniart-gl_2} for a proof). Now, the above proof
shows that irreducible subrepresentations of $U_{1,m}(\chi)$ are
atypical, for $m\geq 1$. This shows the theorem in the present case.
\end{proof}

The pair $(J_{\chi}(1), \chi)$ is a Bushnell-Kutzko type for the
inertial class $s$ (see \cite[Section
8]{Bushnell-kutzko-Semisimpletypes}).  From the above theorem we
deduce the following result:
\begin{theorem}\label{princ_final_rev_cor_thr}
  Let $q_F>3$ if $n>3$ and $q_F>2$ if $n\in \{2,3\}$. Let $\tau$ be a
  typical representation for the inertial class $s=[T_n, \chi]$ then
  $\tau$ is a subrepresentation of
  $\ind_{J_\chi(1)}^{K_n}(\chi)$. Moreover we have
 $$\dim_{\mathbb{C}}\ho_{K_n}(\tau,
 i_{B_n}^{G_n}(\chi))= \dim_{\mathbb{C}}\ho_{K_n} (\tau,
 \ind_{J_{\chi}(1)}^{K_n}(\chi))$$
\end{theorem}
\begin{remark}
  When $\#k_F=2$ and $n=2$ Henniart showed in
  \cite{Henniart-gl_2}[A.2.6, A.2.7] that the Bushnell-Kutzko type for
  the inertial class $s=[T_2, \chi_1\boxtimes \chi_2]$,
  $\chi_1\chi_2^{-1}\neq \id$ has two typical representations one
  given by
$$\ind_{J_{\chi}(1)}^{\g{2}{\integers{F}}}(\chi)$$
and the other representation turns out to be the complement (it
follows from \cite[Proposition 1(b)]{Casselmanres} that
there is a unique complement) of
$\ind_{J_{\chi}(1)}^{\g{2}{\integers{F}}}(\chi)$ in
$\ind_{J_{\chi}(2)}^{\g{2}{\integers{F}}}(\chi)$.  But for $\#k_F>2$
and $n>3$ we expect that typical representations are precisely the
irreducible sub representations of
$$\ind_{J_{\chi}(1)}^{K_n}(\chi).$$
For $\#k_F=2$ and $n>2$ a typical representation may not be contained
in the above representation as shown by Henniart for the case of
$\g{2}{F}$.
\end{remark}
\section{Typical representations 
for \texorpdfstring{$\g{3}{F}$}{}}
Any inertial class of the group $G_3$ belongs to one of the following
classes (see Lemma \ref{standard_form_rev}):
\begin{enumerate}
\item[1] $[G_3, \sigma]$, where $\sigma$ is a cuspidal representation
  of $G_3$
\item[2] $[G_2\times G_1, \sigma\boxtimes\chi]$, where $\sigma$ is a cuspidal
  representation of $G_2$ and $\chi$ is any character of $F^\times$. 
\item[3] $[T_3, \chi=\chi_1\boxtimes\chi_2\boxtimes\chi_3]$, where $\chi_1$,
  $\chi_2$ and $\chi_3$ are three characters of $F^\times$. 
\end{enumerate}
Typical representations for any inertial class of the form
$s=[G_3, \sigma]$ are classified in the work of Pa\v{s}k\={u}nas
\cite{Paskunas-uniqueness}. Up to isomorphism there exists a unique
typical representation for $s$. Similarly, the theorem
\ref{rev_n+1_final} shows that, if $q_F>2$, then up to isomorphism
there exists a unique typical representation for any inertial class of
the form $[G_2\times G_1, \sigma\boxtimes\chi]$. If $q_F>2$, then for
any inertial class $s=[T_3, \chi]$ we showed that any typical
representation occurs as a sub representation of
$\ind_{J_\chi(1)}^K\chi$. The pair $(J_\chi(1), \chi)$ is a
Bushnell--Kutzko type for $s$. Moreover, we also have the multiplicity
result Theorem \ref{princ_final_rev_cor_thr}. We conclude the
following result for $\g{3}{F}$.
\begin{theorem}
  Let $q_F>2$ and $s$ be any inertial class of $G_3$. Typical
  representations for $s$ are precisely the irreducible
  subrepresentations of $\ind_{J_s}^{K_3}\lambda_s$, where
  $(J_s, \lambda_s)$ is a Bushnell--Kutzko type for $s$.
  \end{theorem}

\bibliography{../biblio}
\bibliographystyle{amsalpha}
\noindent
Santosh Nadimpalli, 
IMAPP--Radboud
Universiteit Nijmegen,
Heyendaalseweg 135, 6525AJ Nijmegen, 
The
Netherlands.
\noindent
Email: \texttt{nvrnsantosh@gmail.com}, \texttt{Santosh.Nadimpalli@ru.nl}.
\end{document}